\documentclass[a4paper,article]{LDES}
\usepackage{amsmath,amsfonts}
\usepackage{array}
\usepackage[caption=false,font=normalsize,labelfont=sf,textfont=sf]{subfig}
\usepackage{textcomp}
\usepackage{stfloats}
\usepackage{url}
\usepackage{verbatim}
\usepackage{graphicx}
\usepackage{nomencl}
\usepackage{etoolbox}
\usepackage{ifthen}
\usepackage{booktabs}
\usepackage[ruled]{algorithm2e}
\usepackage{appendix}
\usepackage{multirow}
\usepackage{bm}
\usepackage{makecell}
\usepackage{algpseudocode}
\allowdisplaybreaks 
\usepackage[subtle,tracking=normal]{savetrees}
\usepackage{float}

\usepackage[numbers,sort&compress]{natbib}
\usepackage{nomencl}
\makenomenclature
\renewcommand{\nomgroup}[1]{%
    \ifthenelse{\equal{#1}{A}}{\item[\textbf{Indices and Sets}]}{%
    \ifthenelse{\equal{#1}{B}}{\item[\textbf{Abbreviations}]}{%
    \ifthenelse{\equal{#1}{C}}{\item[\textbf{Parameters}]}{%
    \ifthenelse{\equal{#1}{D}}{\item[\textbf{Decision Variables}]}{}}}}}
\usepackage{multicol}
\usepackage{framed} 
\usepackage{etoolbox} 
\usepackage{float}
\usepackage{placeins}
\usepackage{threeparttable} 

\def\tsc#1{\csdef{#1}{\textsc{\lowercase{#1}}\xspace}}
\tsc{WGM}
\tsc{QE}

\newtheorem{theorem}{Theorem}
\newtheorem{remark}{Remark}

\begin{document}
\let\WriteBookmarks\relax
\def\floatpagepagefraction{1}
\def\textpagefraction{.001}
\let\printorcid\relax 

\shorttitle{\rmfamily N.Qi et al. Long-Term Energy Management of Microgrid}    

\shortauthors{\rmfamily N. Qi et al.}

\title[mode = title]{Long-Term Energy Management for Microgrid with Hybrid Hydrogen-Battery Energy Storage: A Prediction-Free Coordinated Optimization Framework}

\author[1]{Ning Qi}
\cormark[1]
\ead{nq21767@columbia.edu} 
\credit{Conceptualization, Modeling, Methodology, Software, Writing}

\author[2]{Kaidi Huang}
\credit{Methodology, Proof}

\author[1]{Zhiyuan Fan}
\credit{Modeling, Revision}


\author[1]{Bolun Xu}
\credit{Supervision, Revision, Funding Support}


\address[1]{Department of Earth and Environmental Engineering, Columbia University, New York, NY 10027, USA}
\address[2]{Department of Electrical Engineering, Tsinghua University, Beijing 100084, China}

\cortext[1]{Corresponding author} 

\begin{abstract}
This paper studies the long-term energy management of a microgrid coordinating hybrid hydrogen-battery energy storage. We develop an approximate semi-empirical hydrogen storage model to accurately capture the power-dependent efficiency of hydrogen storage. We introduce a prediction-free two-stage coordinated optimization framework, which generates the annual state-of-charge (SoC) reference for hydrogen storage offline. During online operation, it updates the SoC reference online using kernel regression and makes operation decisions based on the proposed adaptive virtual-queue-based online convex optimization (OCO) algorithm. We innovatively incorporate penalty terms for long-term pattern tracking and expert-tracking for step size updates. We provide theoretical proof to show that the proposed OCO algorithm achieves a sublinear bound of dynamic regret without using prediction information. Numerical studies based on the Elia and North China datasets show that the proposed framework significantly outperforms the existing online optimization approaches by reducing the operational costs and loss of load by around 30\% and 80\%, respectively. These benefits can be further enhanced with optimized settings for the penalty coefficient and step size of OCO, as well as more historical references.

\end{abstract}


\begin{highlights}
\item Long-term energy management of microgrid considering seasonal uncertainties and seasonal storage
\item A prediction-free two-stage coordinated optimization framework
\item SoC reference of hydrogen storage generated from kernel regression and historical and AI-generated scenarios
\item A virtual-queue-based online convex optimization algorithm with expert-tracking
\item Numerical studies on Elia and North China with ground-truth datasets spanning 10 years
\end{highlights}

\maketitle

\section{Introduction}

\subsection{Background and motivation}

A microgrid is a self-contained electrical network with resources including energy storage (ES), renewable energy sources (RES), and controllable loads, operated in either grid-connected or island mode~\cite{dey2023microgrid,isolated}. Microgrids enhance energy resilience, promote decarbonization, and reduce transmission system investments, but the volatility of RES poses challenges to short-term supply-demand balances~\cite{pang2023microgrid,microgridreview}. Besides, seasonal variations in RES availability~\cite{guo2023long} and extreme weather events~\cite{zhou2023novel} have highlighted the significance of the long-term energy management of microgrids. 

Hybrid energy storage system (HESS)~\cite{hajiaghasi2019hybrid,qi2023chance} offers a promising way to guarantee both the short-term and long-term supply-demand balance of microgrids.  HESS is composed of two or more ES units with different but complementing characteristics, such as duration and efficiency. 
In day-ahead or intra-day operations, batteries can effectively address the uncertainties introduced by RES and demand. For long-term operation, hydrogen storage consisting of electrolyzer and fuel cell can provide efficient solutions to seasonal energy shifting~\cite{jansen2021cost}. In this paper, we focus on a typical application: hybrid hydrogen-battery energy storage (H-BES). Given the differences in storage properties and unanticipated seasonal uncertainties, designing an effective long-term energy management framework for microgrids with H-BES is significant but challenging.

\subsection{Literature review}

Previous research mainly focuses on the short-term energy management of microgrids with H-BES. Two-stage robust optimization is proposed in~\cite{fan2021robustly} for the market operation of H-BES, where the uncertainties from RES are modeled by uncertainty sets. A two-stage distributionally robust optimization-based coordinated scheduling of an integrated energy system with H-BES is introduced in~\cite{qiu2023two}, where an ambiguity set is employed to model the uncertainties from RES and integrated energy loads. Two-stage stochastic energy management of H-BES is proposed in~\cite{eghbali2022stochastic}, where the uncertainties from RES, load, and prices are modeled by typical scenarios. However, these works rely solely on offline optimization methods with predefined uncertainty modeling, which may face optimality or feasibility issues in real-time operation. This motivates the research on real-time energy management with online optimization methods, such as the rolling-horizon method. Model predictive control (MPC) is the widely used rolling-horizon method and multi-level MPC controllers are developed for microgrids with hydrogen or H-BES in~\cite{trifkovic2013modeling,guo2023long}. An actor-critic deep reinforcement learning method is proposed in~\cite{hu2022soft} to address multi-timescale coordinated dispatch of microgrid with hybrid battery and supercapacitor. MPC and approximate dynamic programming approach are jointly utilized for multi-stage coordinated dispatch~\cite{li2021multi}, which achieves robust real-time performance through continuously updated forecasts. However, the limitations in the aforementioned works mainly lie in \textit{(i)} The short-term energy management methods may face infeasibility issues in the long-term operation when considering seasonal variations RES and load.  \textit{(ii)} The performance of these techniques strongly depends on the accuracy of the prediction of uncertainties. However, the predictions are practically
\nomenclature[A]{$i\text{, }p$}{Indices for parallel learning sequence and hydrogen storage efficiency segment, respectively}
\nomenclature[A]{$s\text{, }t$}{Indices for historical scenarios and time period, respectively}
\nomenclature[A]{$\bm{\Omega}_{{I}}\text{, }\bm{\Omega}_{{P}}$}{Sets for parallel learning sequence and hydrogen storage efficiency segment, respectively}
\nomenclature[A]{$\bm{\Omega}_{{S}}\text{, }\bm{\Omega}_{{T}}$}{Sets for historical scenarios and time period, respectively}
\nomenclature[A]{$\bm{x}\text{, }\bm{\xi}$}{Sets for decision variables and stochastic parameters, respectively}
\nomenclature[B]{ES}{Energy storage}
\nomenclature[B]{HESS}{Hybrid energy storage system}
\nomenclature[B]{H-BES}{Hybrid hydrogen-battery energy storage}
\nomenclature[B]{MPC}{Model predictive control}
\nomenclature[B]{OCO}{Online convex optimization}
\nomenclature[B]{RES}{Renewable energy sources}
\nomenclature[B]{RMSE}{Root mean square error}
\nomenclature[B]{Reg}{Regret for OCO algorithm}
\nomenclature[B]{SDP}{Stochastic dynamic programming}
\nomenclature[B]{SoC}{State of charge}
\nomenclature[C]{$A_p \text{, } B_p$}{Slope and the intercept of piecewise linear charging segments, respectively}
\nomenclature[C]{$C_p \text{, } D_p$}{Slope and the intercept of piecewise linear discharging segments, respectively}
\nomenclature[C]{$\overline{P}^{\text{B}}\text{, }\overline{P}^{\text{H}}$}{Upper power bounds of battery and hydrogen storage, respectively}
\nomenclature[C]{$\underline{P}^{\text{B}}$}{Lower power bound of battery}
\nomenclature[C]{$\overline{E}^{\text{B}}\text{, }\overline{E}^{\text{H}}$}{Upper SoC bounds of battery and hydrogen storage, respectively}
\nomenclature[C]{$\underline{E}^{\text{B}}\text{, }\underline{E}^{\text{H}}$}{Lower SoC bounds of battery and hydrogen storage, respectively}
\nomenclature[C]{$\eta^{\text{B,c}}$\text{, }$\eta^{\text{B,d}}$}{Charging and discharging efficiency of battery storage, respectively}
\nomenclature[C]{$\eta^{\text{H,c}}$\text{, }$\eta^{\text{H,d}}$}{Charging and discharging efficiency of hydrogen storage, respectively}
\nomenclature[C]{$\varepsilon$}{Self-discharge rate of battery storage}
\nomenclature[C]{$\underline{P}^{\text{D}}\text{, }\overline{P}^{\text{D}}$}{Lower and upper power bounds of diesel generator, respectively}
\nomenclature[C]{$RD^{\text{D}}\text{, }RU^{\text{D}}$}{Downward and upward ramping rates of  diesel generator, respectively}
\nomenclature[C]{$c^\text{B}\text{, }c^\text{H}\text{, }c^\text{D}\text{, }c^\text{L}$}{Marginal discharge costs of battery and hydrogen storage, fuel price diesel generator, load curtailment price, respectively}
\nomenclature[C]{$\omega_{s,t}\text{, }\varphi$}{Dynamic weights for historical scenarios and penalty coefficient for reference tracking, respectively}
\nomenclature[C]{$\alpha_{i,t}\text{, }\beta_{i,t}\text{, }\gamma$}{step sizes for OCO algorithm, respectively}
\nomenclature[C]{$Q_{i,t}$}{Virtual queue for OCO algorithm}
\nomenclature[D]{$E_{t}^{\text{B}}\text{, }E_{t}^{\text{H}}$}{State of charge of battery and hydrogen storage, respectively}
\nomenclature[D]{$P_{t}^{\text{B,c}}\text{, }P_{t}^{\text{B,d}}$}{Charging and discharging power of battery storage, respectively}
\nomenclature[D]{$P_{t}^{\text{H,c}}\text{, }P_{t}^{\text{H,d}}$}{Charging and discharging power of hydrogen storage, respectively}
\nomenclature[D]{$z_{t}^{\text{c}}\text{, }z_{t}^{\text{d}}$}{Binary variables for charging and discharging segments of hydrogen storage, respectively}
\nomenclature[D]{$h_{t}^{\text{c}}\text{, }h_{t}^{\text{d}}$}{Hydrogen production and consumption, respectively}
\nomenclature[D]{$P_{t}^{\text{D}}\text{, }P_{t}^{\text{L}}\text{, }P_{t}^{\text{R}}$}{Diesel generator power, loss of load power and dispatched renewable power, respectively}
\twocolumn[{%
    \begin{center}
    \end{center}
    \begin{framed}
    \begin{multicols}{2}
    \printnomenclature
    \end{multicols}
    \end{framed}
}]

\noindent unavailable or unreliable for microgrid operators.

To address the first limitation, recent studies have started to explore the long-term energy management of microgrids, which aims to solve the multi-time-period dispatch with non-anticipativity. Stochastic dynamic programming is technically sound, which can decompose the multi-period dispatch problem into sequential single-period dispatch problems through value function. And it is applied in~\cite{darivianakis2017data} by learning the value function of H-BES. However, it becomes computationally intractable to train the value function if the storage duration spans multiple months. A continuous spectrum splitting approach is proposed in~\cite{feng2024hybrid} to assign low-frequency uncertainty scenarios to hydrogen and high-frequency uncertainty scenarios to batteries for power balance, but this approach is designed for the planning of H-BES. A data-driven coordinated dispatch framework is proposed in~\cite{guo2023long}, where the state of charge (SoC) reference for hydrogen storage is generated based on historical simulations. This reference is then updated and embedded into MPC for real-time operation. However, the use of MPC makes the entire framework dependent on forecasting. 

Additionally, prediction-free online optimization methods are gaining increased attention. Lyapunov optimization
and online convex optimization (OCO) are effective representatives~\cite{wang2023online}. Lyapunov optimization adopts a ``1-lookahead'' pattern, where uncertainties are observed first, followed by solving the Lyapunov drift problem~\cite{shi2015real}. It has wide applications in demand response~\cite{zheng2014distributed}, electric vehicle charging~\cite{yan2023real}, microgrid~\cite{alzahrani2023real}, etc. The long-term operational cost minimization of hydrogen-based building energy systems is transformed into several single-slot subproblems using Lyapunov optimization~\cite{yu2022joint}. A joint energy scheduling and trading algorithm based on Lyapunov optimization and a double-auction mechanism is designed in~\cite{zhu2020energy} to optimize the long-term energy cost of each microgrid. However, in some cases, the uncertainties can not be observed before decision-making and 
Lyapunov optimization becomes inapplicable. For instance, storage participants bid with unknown future prices, and the prices are cleared by the market after the bidding process~\cite{zheng2023energy}.
Instead, OCO adopts a ``0-lookahead'' pattern, where the decision is made before the observation of uncertainties. 
And OCO has been utilized in demand side management~\cite{kim2016online} and ancillary services~\cite{zhao2020distributed} due to its completely prediction-free and fast response nature. However, to the best of our knowledge, no research has addressed the long-term energy management of microgrids with H-BES within the OCO framework. The application of OCO in the focused topic may face the following challenges: \textit{(i)} OCO is problem-dependent without a predefined mathematical formulation, and there is no prior experience available as a reference for designing OCO for microgrids with H-BES. \textit{(ii)} Although recent works~\cite{muthirayan2022online,liu2022simultaneously,yi2022regret,ding2021dynamic} have embedded inter-temporal constraints into the OCO framework, OCO still risks falling into local optima due to its myopic nature. \textit{(iii)} OCO aims to achieve regret (Reg) that grows sublinearly with time horizon $T$. However, most of the existing OCO algorithms fail to address the sublinear bounds for dynamic Reg~\cite{liu2022simultaneously,yi2022regret,ding2021dynamic} or require prediction information to improve the performance~\cite{muthirayan2022online}. Please see Table~\ref{literature review} for a comprehensive comparison.

\subsection{Research gap}

Existing literature is summarized in Table~\ref{literature review}. Although some works achieve good results in the long-term energy management of microgrids with H-BES, there are still several research gaps that have not been adequately addressed.

(1) Most existing studies employ a simplified operational model for hydrogen storage, using a constant energy conversion efficiency regardless of whether the storage operates at full capacity. However, the efficiency of hydrogen storage varies with the charge/discharge power and follows a nonlinear function~\cite{ulleberg2003modeling}. Using a simplified model can result in sub-optimal or even infeasible solutions~\cite{baumhof2023optimization}. Therefore, it is crucial to incorporate this nonlinearity into the microgrid energy management with H-BES.

(2) Current microgrid energy management approaches either employ offline optimization methods (e.g., robust optimization~\cite{fan2021robustly}, frequency-domain method~\cite{feng2024hybrid}) or prediction-dependent online optimization methods (e.g., MPC~\cite{guo2023long}, stochastic dynamic programming~\cite{darivianakis2017data}). However, the distribution and prediction information is often inaccurate or unavailable in practical microgrid operations. Thus, designing a prediction-free optimization framework for microgrid energy management with H-BES is necessary.

(3) OCO is a promising ``0-lookahead'' online optimization method originating from the fields of machine learning and control~\cite{yi2022regret}-\cite{ding2021dynamic}. However, OCO lacks a global view of long-term operations and adaptability to the high volatility of microgrids. Hence, it is important to extend traditional OCO methods to incorporate long-term operational patterns and time-varying properties. 
\begin{table*}[!ht]
\footnotesize\rmfamily
  \centering
  \begin{threeparttable}
  \caption{\rmfamily Comparison of existing literature on long-term and short-term energy management of H-BES.}
  \setlength{\tabcolsep}{0.4mm}{
      \begin{tabular}{c c c c c}
    \toprule
    Reference  & Storage Type \& Model & Long-term Optimization & Short-term Optimization & Prediction-Free \\
    \midrule
    \cite{fan2021robustly}  & H-BES-Constant & X     & Robust Optimization &  $\checkmark$(Offline) \\
    \cite{qiu2023two}  & H-BES-Constant & X     & Distributionally Robust Optimization & $\checkmark$(Offline) \\
    \cite{eghbali2022stochastic}  & H-BES-Constant & X     & Stochastic Optimization & $\checkmark$(Offline) \\
    \cite{trifkovic2013modeling}  & H-BES-Electrochemical & X     & MPC & X \\
    \cite{hu2022soft}  & Battery+Supercapacitor-Constant & X     & Deep Reinforcement Learning & X \\
    \cite{li2021multi}  & Battery+Thermal Storage-Constant & X     & MPC+Dynamic Programming & X \\
    \cite{darivianakis2017data}  & H-BES-Constant & \multicolumn{2}{c}{Stochastic Dynamic Programming} & X \\
    \cite{feng2024hybrid}  & H-BES-Constant & \multicolumn{2}{c}{Spectrum Splitting Approach} & $\checkmark$(Offline) \\
    \cite{guo2023long}   & H-BES-Constant & Historical Reference & MPC & X \\
    \cite{yu2022joint}  & H-BES-Constant & \multicolumn{2}{c}{Lyapunov Optimization} & $\checkmark$(1-lookahead) \\
    \cite{zhu2020energy} & Hydrogen Full Cell-Constant & \multicolumn{2}{c}{Lyapunov Optimization} & $\checkmark$(1-lookahead) \\
        \cite{muthirayan2022online}   & Not Given & X &OCO: Dynamic Reg $\mathcal{O}(T^{\max \{1-a-c, c\}})$, 0<a,c<1 & $\checkmark$(0-lookahead) \\    \cite{liu2022simultaneously}  & Not Given & X & OCO: Dynamic Reg $\mathcal{O}(T^{c} P_x^{c})$, 0<c<1 & $\checkmark$(0-lookahead) \\
    \cite{yi2022regret} & Not Given & X &OCO: Static Reg $\mathcal{O}(T^{\max \{1-c, c\}})$, 0<c<1 & $\checkmark$(0-lookahead) \\
    \cite{ding2021dynamic}  & Not Given & X &OCO: Dynamic Reg $\mathcal{O}(\max (T^c P_x, T^{1-c}))$, 0<c<1 & $\checkmark$(0-lookahead) \\
    This Paper & H-BES-Semi-Empirical & Historical\&AI-Generated Reference & \makecell{OCO: Dynamic Reg $\mathcal{O}(T^c (1 + P_x)^{1 - \kappa} + T^{1 - c} (1 + P_x)^{\kappa})$,\\ $0<\kappa<c<1$} & $\checkmark$(0-lookahead) \\
    \bottomrule
    \end{tabular}%
    }\label{literature review}
    \begin{tablenotes}
\item[a] Depending on the baselines used, Reg is divided into static Reg, with the baseline being a single-period optimal solution, and dynamic Reg, with the baseline being the global optimal solution.
\item[b] $P_{x}$: path-length, i.e., the accumulated variation of optimal decisions; $P_{g}$: function variation, i.e., the accumulated variation of constraints.
\end{tablenotes}
\end{threeparttable}\vspace{-0.5cm}
\end{table*}%

\subsection{Contributions}

Motivated by the research gaps, this paper proposes a prediction-free coordinated optimization framework for long-term energy management of microgrid with H-BES while incorporating the nonlinearity of hydrogen storage and seasonal uncertainties from RES and load. Specifically, our contributions are threefold:

\textbf{(1) Modeling:} We propose an approximate semi-empirical hydrogen storage model using piecewise linear relaxation, which accurately captures the power-dependent efficiency of hydrogen storage. Simulations demonstrate that, compared to the constant efficiency model, the proposed approximation model avoids both overly optimistic and overly conservative strategies. This results in a reduction of the practical yearly operational cost by 10\% or 36\%, and a decrease in yearly loss of load by 1.94 MWh or 3.85 MWh.

\textbf{(2) Solution Methodology:} We introduce a prediction-free two-stage coordinated optimization framework. In the offline stage, the ex-post SoC references for hydrogen storage are generated by deterministic mixed-integer linear programming with historical and AI-generated data on RES and load. These references help to avoid myopic online decision-making and are incrementally updated by kernel regression with newly observed data. Subsequently, we develop an adaptive virtual-queue-based OCO algorithm for prediction-free online decision-making. Compared to the traditional OCO algorithm~\cite{muthirayan2022online,liu2022simultaneously,yi2022regret,ding2021dynamic}, the proposed method innovatively incorporates a penalty term for long-term pattern tracking and expert-tracking for step size updates. The proposed OCO algorithm is proven to achieve a sublinear bound for dynamic regret.

\textbf{(3) Numerical Study:} We demonstrate the effectiveness of the proposed framework using ground-truth data from Elia~\cite{Elia-data} and North China~\cite{North-China-data}. Simulations show that introducing the reference significantly reduces operational costs and loss of load by 40\%-57\% and 60\%-90\%, respectively. Furthermore, compared to the prediction-dependent MPC method, the prediction-free OCO method further decreases operational costs and loss of load by 24\%-29\% and 73\%-89\%, respectively. These benefits can be further enhanced with optimized settings for the penalty coefficient and step size of OCO, as well as more historical references.

\subsection{Paper Organization}

We organize the remainder of the paper as follows. Section~\ref{hydrogenmodel} presents an approximate semi-empirical modeling of hydrogen storage. Section~\ref{PF} provides the problem formulation for long-term energy management of the microgrid with H-BES. Section~\ref{framework} introduces the prediction-free two-stage coordinated optimization framework and the proof of OCO performance. Section~\ref{case} describes numerical case studies to verify the effectiveness of the proposed framework. Finally, we conclude this paper in Section~\ref{conclusion}.

\section{Approximate semi-empirical hydrogen energy\\ storage model}~\label{hydrogenmodel}
\subsection{Structure of hydrogen storage system}

\begin{figure}[!ht]
\vspace{-0.5em}
  \footnotesize\rmfamily   \setlength{\abovecaptionskip}{-0.1cm}  
    \setlength{\belowcaptionskip}{-0.1cm} 
  \begin{center}  \includegraphics[width=0.9\columnwidth]{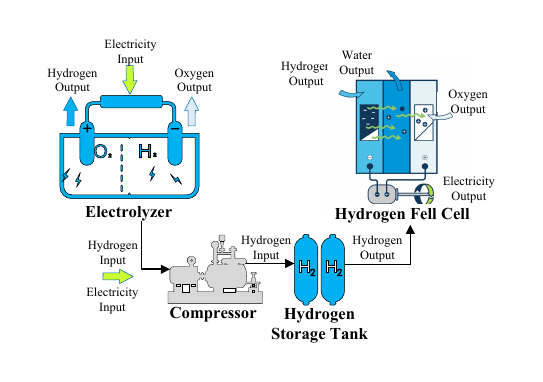}
     \caption{\rmfamily Schematic diagram of hydrogen storage system.}\label{HESS}
  \end{center}
  \vspace{-0.5cm}
\end{figure}

A hydrogen storage system is composed of several key components, such as electrolyzers, hydrogen storage tanks, fuel cells, compressors, and other auxiliary equipment, as illustrated in Fig.~\ref{HESS}. Electrolyzers convert electrical energy into chemical energy by producing hydrogen and oxygen. This paper considers the most mature and commonly used alkaline water electrolyzer. Hydrogen storage tanks are used to store the produced hydrogen. Fuel cells convert the stored hydrogen back into electricity, and we consider the typical type, proton exchange membrane fuel cell (PEMFC). Other auxiliary equipment, including the compressor, cooling system, and control system, is excluded from the modeling.

\subsection{Alkaline water electrolyzer model}

(1) Polarization curve

The polarization curve describes the electrochemical behavior of an electrolyzer, modeling the relationship between current and voltage. To account for the impact of temperature and pressure on the thermodynamics and electrochemical process within the electrolyzer, we combine the most used model proposed by Ulleberg~\cite{ulleberg2003modeling} and the modified model proposed by Sanchez~\cite{sanchez2018semi}:
\begin{equation}
      \begin{split}
    U_{\text{cell}}^{\text{E}} &= U_{\mathrm{rev}}  + \left[\left(r_1 + d_1\right) + r_2 \cdot \theta + d_2 \cdot P\right] \cdot \dfrac{i}{A} \\
    & + s \cdot \log\left[\left(t_1 + \frac{t_2}{\theta} + \frac{t_3}{\theta^2}\right) \cdot \dfrac{i}{A} + 1\right] 
    \end{split}\label{U-I}
\end{equation}
where the reversible voltage and cell voltage of the electrolyzer are defined as~$U_{\mathrm{rev}}$ and~$U_{\text{cell}}^{\text{E}}$. Temperature and pressure are given by $\theta$ and $P$. The current and effective area of the electrode is defined as~$i$ and $A$. Parameters $r_1$, $r_2$, $d_1$, $d_2$, $t_1$, $t_2$, $t_3$, $s$ are the constants which can be learned from the experimental data.

(2) Faraday efficiency

Faraday efficiency is defined as the ratio of measured hydrogen production to the theoretical value. For an alkaline electrolyzer, the Faraday efficiency typically ranges from 85\% to 95\% and is affected by temperature. We adopt the four-parameter Faraday efficiency model as~\eqref{Faraday}. 
\begin{equation}
\eta_{\text{F}}=\left(\frac{(i/A)^2}{f_1+f_2 \cdot \theta+(i/A)^2}\right) \cdot\left(f_3+f_4 \cdot \theta\right)\label{Faraday}
\end{equation}
where Faraday efficiency is defined as~$\eta_{\text{F}}$. Parameters $f_1$, $f_2$, $f_3$, $f_4$ are the constants which can be learned from the experimental data.

(3) Approximate charging efficiency

According to Faraday's law, the hydrogen production rate is defined as~\eqref{h2pro}. The charging efficiency is given by~\eqref{h2charge}.
\begin{subequations}
    \begin{align}
      &h^{\text{c}}=3600 \cdot \dfrac{\eta_{\text{F}}\cdot M \cdot i\cdot N}{2 F}\label{h2pro}\\
      & \eta^{\text{H},{\text{c}}}=\dfrac{h^{\text{c}} \cdot \text{LHV}}{P_{\text{Stack}}}=3600 \cdot \dfrac{\eta_{\text{F}}\cdot  M\cdot  \text{LHV}}{2 F \cdot U_{\text{cell}}}\label{h2charge}
    \end{align}
\end{subequations}
where $h^{\text{c}}$ is the hydrogen production rate of electrolyzer. $M$ is the molar mass of hydrogen. $F$ is the Faraday's constant, i.e., 96485 C/mol. $N$ is the number of cells of the stack. $\text{LHV}$ is the lower heat value of hydrogen, i.e., 33.33 kWh/kg.

\begin{figure}[!ht]
\vspace{-0.5em}
 \footnotesize\rmfamily     \setlength{\abovecaptionskip}{-0.1cm}  
    \setlength{\belowcaptionskip}{-0.1cm} 
  \begin{center}  \includegraphics[width=1\columnwidth]{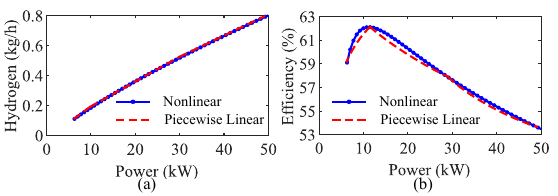}
     \caption{\rmfamily Approximation of alkaline electrolyzer properties at $90\,^\circ\mathrm{C}$, 10 bar: (a) hydrogen production, (b) efficiency.}\label{chargeapproximate}
  \end{center}
  \vspace{-0.5cm}
\end{figure}

As illustrated in Fig.~\ref{chargeapproximate}, the blue curves from the semi-empirical model are non-linear and power-dependent, including a peak in efficiency at around 20\% of the rated power. Therefore, a constant conversion efficiency cannot capture the variations in efficiency. To facilitate dispatch optimization, we adopt a piecewise linear approximation for hydrogen production, depicted by red dashed lines.

\subsection{PEMFC model}

(1) Polarization curve

The polarization curve of PEMFC is typically modeled using the equivalent circuit model proposed by Amphlett~\cite{amphlett1995performance}. The cell voltage $U_{\text {cell}}^{\text{F}}$ is given by~\eqref{PEMFC}, which equals the open circuit voltage $E_{\text {Nernst }}$ dropped by three types of irreversible losses: activation losses $U_{\text {act }}$, ohmic losses $U_{\text {ohmic }}$, and concentration losses $U_{\text {con }}$.
\begin{subequations}
    \begin{align}
    &  U_{\text {cell}}^{\text{F}}=E_{\text {Nernst }}-U_{\text {act }}-U_{\text {ohmic }}-U_{\text {con }}\label{vcell}\\
     & E_{\text {Nernst }}=\frac{1}{2 F}\big[\Delta G-\Delta S (\theta-\theta_{\text {ref }})\nonumber\\
     &\hspace{1cm}+R \cdot \theta\left(\log( P_{\mathrm{H}_2})+\frac{\log (P_{\mathrm{O}_2})}{2}\right)\big]\label{Enernst}\\
    & U_{\text {act }}=a_1+a_2 \cdot \theta+a_3 \cdot \theta \cdot \log (C_{\mathrm{O}_2})+a_4 \cdot \theta \cdot \log \left(i\right)\label{Vact}\\
        & U_{\text {ohmic }}=i \cdot R_{\text {ohmic }}=i\left(r_M\cdot l/A+R_c\right)\label{Vohm}\\
            & U_{\text {con }}=B \cdot \log \left(1-\dfrac{J}{J_{\max }}\right)\text{, }J=\dfrac{i}{A}\label{VCON}
    \end{align}\label{PEMFC}
\end{subequations}
where $\Delta G$ is the Gibbs free energy. $\Delta S$ is the entropy change. $R$ is the gas constant (8.314 J/(K$\cdot$mol)). $P_{\mathrm{H}_2}$ and $P_{\mathrm{O}_2}$ are the partial pressures of hydrogen and oxygen respectively. $T_{\text{ref}}$ is the reference temperature (298.15 K).  $C_{\mathrm{O}_2}$ is the oxygen concentration at the surface of the cathode catalyst. $r_M$ is the resistivity of the electrolyte membrane. $l$ is the thickness of the electrolyte membrane. $B$ is the concentration overpotential coefficient. $J$ and $J_{\text{max}}$ are the current density and its maximum value. $a_1$, $a_2$, $a_3$, $a_4$ are constants that can be learned from the experimental data.

(2) Approximate discharging efficiency

According to Faraday's law, the hydrogen consumption rate is defined as~\eqref{h2cos}. The discharging efficiency is given by~\eqref{h2discharge}.
\begin{subequations}
    \begin{align}
      &h^{\text{d}}=3600 \cdot \dfrac{ M \cdot i\cdot N}{2 F}\label{h2cos}\\
      & \eta^{\text{H},{\text{d}}}=\dfrac{P_{\text{Stack}}}{h^{\text{d}} \cdot \text{HHV}}= \dfrac{2 F \cdot U_{\text{cell}}}{3600 M\cdot  \text{HHV}}\label{h2discharge}
    \end{align}
\end{subequations}
where $h^{\text{d}}$ is the hydrogen consumption rate of PEMFC. $\text{HHV}$ is the higher heat value of hydrogen, i.e., 39.4 kWh/kg.

\begin{figure}[!ht]
\vspace{-0.5em}
  \footnotesize\rmfamily    \setlength{\abovecaptionskip}{-0.1cm}  
    \setlength{\belowcaptionskip}{-0.1cm} 
  \begin{center}  \includegraphics[width=1\columnwidth]{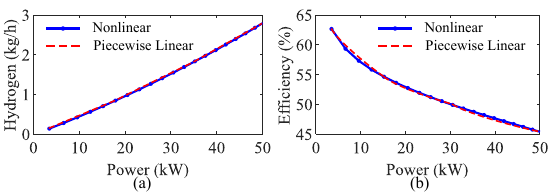}
     \caption{\rmfamily Approximation of PEMFC properties at $60\,^\circ\mathrm{C}$, 10 bar: (a) hydrogen consumption, (b) efficiency.}\label{dischargeapproximate}
  \end{center}
  \vspace{-0.5cm}
\end{figure}

The blue curves from the semi-empirical model in Fig.~\ref{dischargeapproximate} are non-linear and power-dependent. We also adopt a piecewise linear approximation for hydrogen consumption, depicted by red dashed lines.

\subsection{Equivalent hydrogen storage model}

The equivalent hydrogen storage model is presented in~\eqref{hydrogen}. Constraint~\eqref{E-power} defines the relationship between SoC, charge power, and discharge power. Constraints~\eqref{E-bound} limit the SoC of hydrogen storage within the bounds. Constraint~\eqref{E-balance} guarantees ensures a sustainable energy state for hydrogen storage over cycles. Constraints~\eqref{piecewise}-\eqref{piece0-1} describe the tractable formulation of piecewise linear charging and discharging functions. Constraints~\eqref{piecepowerbound} limit hydrogen storage's charging and discharging power. 

\noindent\textbf{ Constraints:} $\forall t\in {{\bm{\Omega} }_{T}}~\forall p\in {{\bm{\Omega} }_{{P}}}$
\begin{subequations}
\begin{align}
  &\hspace{-0.5cm} E_{t+1}^{\text{H}}=E_{t}^{\text{H}}+\Delta t(h_{t}^{\text{c}}-h_{t}^{\text{d}})-E_{t}^{\text{H,L}} \label{E-power}\\ 
 & \hspace{-0.5cm} \underline{E}^{\text{H}}\le E_{t}^{\text{H}}\le \overline{E}^{\text{H}}\label{E-bound}\\
 & \hspace{-0.5cm} E_{T}^{\text{H}}\geq E_{0}^{\text{H}}\label{E-balance}\\
  & \hspace{-0.5cm} h_{t}^{\text{c}}=\sum_p (A_p P_{p\text{,}t}^{\text{H\text{,}c}}+B_p z_{p\text{,}t}^{\text{c}})\text{, }h_{t}^{\text{d}}=\sum_p (C_p P_{p\text{,}t}^{\text{H,d}}+D_p z_{p\text{,}t}^{\text{d}})\label{piecewise}\\
   & \hspace{-0.5cm} P_{t}^{\text{H\text{,}c}}=\sum\nolimits_p P_{p\text{,}t}^{\text{H\text{,}c}}\text{, }P_{t}^{\text{H,d}}=\sum\nolimits_p P_{p\text{,}t}^{\text{H,d}}\label{combinedpower}\\ 
    & \hspace{-0.5cm} \sum\nolimits_p z_{p\text{,}t}^{\text{c}}=1\text{, }\sum\nolimits_p z_{p\text{,}t}^{\text{d}}=1\label{piece0-1}\\
        & \hspace{-0.5cm} \underline{P}_{p}^{\text{H}}z_{p\text{,}t}^{\text{c}}\le P_{p\text{,}t}^{\text{H\text{,}c}}\le \overline{P}_{p}^{\text{H}}z_{p\text{,}t}^{\text{c}}\text{, }0\le P_{p\text{,}t}^{\text{H,d}}\le \overline{P}_{p}^{\text{H}}z_{p\text{,}t}^{\text{d}}\label{piecepowerbound}
\end{align}\label{hydrogen}\vspace{-0.3cm}
\end{subequations}

\noindent where ${{\bm{\Omega} }_{T}}$ and ${{\bm{\Omega} }_{S}}$ are the set of time and parameter segments, respectively. $P_{t}^{\text{H,c}}$, $P_{t}^{\text{H,d}}$, and $E_{t}^{\text{H}}$ are decision variables for the charge power, discharge power, and SoC of hydrogen storage. The SoC of hydrogen storage can be measured by the hydrogen mass or as a ratio of the rated capacity. $E_{t}^{\text{H,L}}$ is the hydrogen load for industrial production processes, such as fertilizer manufacturing and steel-making. $\underline{E}^{\text{H}}$ and $\overline{E}^{\text{H}}$ are the lower and upper bounds of SoC. $\underline{P}^{\text{H}}$ and $\overline{P}^{\text{H}}$ are the lower and upper bounds of power. The lower charging power bound is set by the minimum operating power of the electrolyzer, typically 15\%-20\% of the nominal power. $A_p$ and $B_p$ are the slope and the intercept of piecewise linear charging segments. $C_p$ and $D_p$ are the slope and the intercept of piecewise linear discharging segments. $z_{p\text{,}t}^{\text{c}}$ and $z_{p\text{,}t}^{\text{d}}$ are binary variables for piecewise linear function.

\section{Long-term energy management of microgrid}~\label{PF}

\subsection{Microgrid structure}

\begin{figure}[!ht]
\vspace{-0.5em}
 \footnotesize\rmfamily     \setlength{\abovecaptionskip}{-0.1cm}  
    \setlength{\belowcaptionskip}{-0.1cm} 
  \begin{center}  \includegraphics[width=0.9\columnwidth]{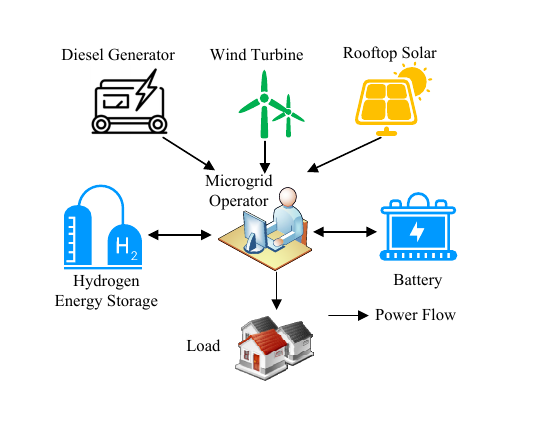}
     \caption{\rmfamily Diagram of microgrid structure.}\label{microgrid}
  \end{center}
  \vspace{-0.8cm}
\end{figure}

In this paper, we only consider the island mode of the microgrid, and the microgrid structure is illustrated in Fig.~\ref{microgrid}. The microgrid consists of renewable generators (wind and solar), diesel generators, H-BES and local loads.

\subsection{Problem formulation}
The objective defined in~\eqref{object} aims to minimize the system cost. This cost comprises the production costs of the diesel generator, penalties for load curtailment (island mode), and operational costs of H-BES. Constraints~\eqref{DGpower} and~\eqref{DGramp} define the power bounds and ramping bounds of the diesel generator. Constraints~\eqref{battery} define the constraints for battery, which are similar in formulation to those for hydrogen storage~\eqref{hydrogen}, as both types of storage involve constraints on charging and discharging rates, SoC, etc. However, it is important to note that the battery efficiency is considered to be constant, there is no minimum charging power limit in~\eqref{power-bound1-B}, and the self-discharge rate should be considered in~\eqref{E-power-B}. Constraints~\eqref{loadRES} limit the load curtailment and dispatchable RES. Constraint~\eqref{feeder} limits the power import from the main grid. Power balance constraint is defined as~\eqref{powerbalance}. The complementary constraints for charging and discharging of battery and hydrogen storage are relaxed and have been removed from the model since sufficient conditions are satisfied~\cite{li2015sufficient}, i.e., discharging price (``+'') is greater than the charging price (``0''). Moreover, the power flow constraints are overlooked within the dispatch model since the microgrid network is generally designed with high reliability and large redundancy~\cite{shuai2018stochastic}.

\noindent\textbf{Objective Function:}
\begin{subequations}\label{object}
\begin{align}
    &\underset{\bm{x}}{\min}\ G(\bm{x}, {\bm{\xi}})=\sum_{t\in\boldsymbol{\Omega}_T}\left(C_t^\text{L}+C_t^\text{D}+C_t^\text{B}+C_t^\text{H}\right)
   \label{costall} \\
&C_t^\text{L}=c^{\text{L}}P_{t}^{\text{L}}\Delta t\text{, }C_t^\text{D}=c^\text{D}P_{t}^\text{D}\Delta t\text{, }C_t^\text{B/H}=c^\text{B/H}P_{t}^\text{B/H,d}\Delta t\label{grid cost}
\end{align}
\end{subequations}
where $c^{\text{L}}$ and $P_{t}^{\text{L}}$ are the load curtailment price and load curtailment power. $c^\text{D}$ and $P_{t}^{\text{D}}$ are the fuel price and power of diesel generator. $c^\text{B}$ and $c^\text{H}$ are marginal discharge costs of battery and hydrogen storage. $\Delta t$ is the time interval. 

\noindent\textbf{Constraints:} $\forall t\in {{\bm{\Omega} }_{T}}~$
\begin{subequations}\label{DG}
\begin{flalign}
& \underline{P}^{\text{D}}\leq P_{t}^{\text{D}}\leq\overline{P}^{\text{D}} && \label{DGpower} \\
& -RD^{\text{D}}\leq P_{t+1}^{\text{D}}-P_{t}^{\text{D}}\leq RU^{\text{D}} && \label{DGramp}
\end{flalign}
\end{subequations}
\vspace{-0.8cm}
\begin{subequations}\label{battery}
\begin{flalign}
& E_{t+1}^{\text{B}}=(1-\varepsilon \Delta t)E_{t}^{\text{B}}+\Delta t(\eta^{\text{B,c}}P_{t}^{\text{B,c}}-P_{t}^{\text{B,d}}/{\eta^{\text{B,d}}}) && \label{E-power-B}\\ 
& \underline{E}^{\text{B}}\le E_{t}^{\text{B}}\le \overline{E}^{\text{B}} && \label{E-bound-B}\\
& E_{T}^{\text{B}}\geq E_{0}^{\text{B}} && \label{E-balance-B}\\
& 0\le P_{t}^{\text{B,c}}\le \overline{P}^{\text{B}} && \label{power-bound1-B}\\
& 0\le P_{t}^{\text{B,d}}\le \overline{P}^{\text{B}} && \label{power-bound2-B}
\end{flalign}
\end{subequations}
\vspace{-0.8cm}
\begin{subequations}\label{loadRES}
\begin{flalign}
& 0\leq P_{t}^{\text{L}}\leq \xi_{t}^{\text{L}} && \label{Loadcurtailment}\\
& 0\leq P_{t}^{\text{R}}\leq \xi_{t}^{\text{R}} && \label{REScurtailment}
\end{flalign}
\end{subequations}
\vspace{-0.8cm}
\begin{flalign}\label{feeder}
& 0\le P_{t}^{\text{G}}\le \overline{P}_{t}^{\text{G}} &&
\end{flalign}
\vspace{-0.8cm}
\begin{flalign}\label{powerbalance}
& P_{t}^{\text{G}}+P_{t}^{\text{R}}+(P_{t}^{\text{B,d}}-P_{t}^{\text{B,c}})+(P_{t}^{\text{H,d}}-P_{t}^{\text{H,c}})+P_{t}^{\text{L}}= \xi_{t}^{\text{L}} &&
\end{flalign}
where $\underline{P}^{\text{D}}$ and $\overline{P}^{\text{D}}$ are the lower and upper power bounds of diesel generator. $RD^{\text{D}}$ and $RU^{\text{D}}$ are downward and upward ramping rates of  diesel generator. $P_{t}^{\text{B,c}}$, $P_{t}^{\text{B,d}}$, and $E_{t}^{\text{B}}$ are decision variables for the charge power, discharge power, and SoC of battery. $\eta^{\text{B,c}}$ and $\eta^{\text{B,d}}$ are the charge and discharge efficiency of battery. $\varepsilon$ is the self-discharge rate of battery. $\underline{E}^{\text{B}}$ and $\overline{E}^{\text{B}}$ are the lower and upper SoC bounds of battery. $\overline{P}^{\text{B}}$ is the upper power bound of battery. $\xi_{t}^{\text{L}}$ and $\xi_{t}^{\text{R}}$ are the load power and available RES power with uncertainties. $P_{t}^{\text{R}}$ is the dispatched RES power. The set of stochastic parameters is given by $\bm{\xi}=\{\xi_{t}^{\text{L}},\xi_{t}^{\text{R}}\}$. The set of decision variables is given by $\bm{x}=\{P_{t}^{\text{L}},P_{t}^{\text{D}},P_{t}^{\text{R}},P_{t}^{\text{B,c/d}},P_{t}^{\text{H,c/d}},E_{t}^{\text{B}},E_{t}^{\text{H}},h_{t}^{\text{c/d}},z_{t}^{\text{c/d}}\}$.

The multi-time-period economic dispatch of microgrid with H-BES ($\textbf{P}_{1}$) is summarized in~\eqref{P1}. Next, we present the methodology for solving this problem. 
\begin{alignat}{1}\label{P1}
  (\textbf{P}_{1})\hspace{4pt} \underset{\bm{x}}{\min} & {\hspace{4pt}} G(\bm{x},{\bm{\xi}})\notag\\ 
 \text{s.t.}  & {\hspace{4pt}} \eqref{hydrogen}\text{, }\eqref{DG}-\eqref{powerbalance}
\end{alignat}
\section{Prediction-free coordinated optimization framework}~\label{framework}

\subsection{Motivations}
Solving the problem~\eqref{P1} has the following challenges:

\textbf{(1) Non-anticipatively :} The long-term energy management of the microgrid typically spans more than one month or one season. Nevertheless, the forecast accuracy is acceptable only for several hours ahead. Hence, the load power and available RES power are unanticipated in the long-term optimization. And online optimization methods should be adopted to decompose the long-term optimization problem into several short-term optimization problems.

\textbf{(2) Storage Dispatch Priority:} Batteries with lower marginal discharge costs will be given priority over hydrogen storage with higher marginal discharge costs. We defer the complete proof to Appendix~\ref{priority}. The battery-prioritized strategy is feasible and economical for short-term operation. However, this approach does not account for seasonal variations in RES and load, which will result in a lack of pre-stored hydrogen and load losses in long-term operations. Therefore, it is necessary to design a ``reference'' with a global view to help guide hydrogen storage actions.

\textbf{(3) Convexity:} The piecewise linearization will introduce nonconvexity to the optimization, which contradicts the overall logic of most convex optimization approaches. However, introducing a global ``reference'' can mitigate this challenge by pre-determining the efficiency.

\subsection{Two-stage coordinated optimization framework}
\begin{figure*}[!ht]
\vspace{-0.5em}
 \footnotesize\rmfamily     \setlength{\abovecaptionskip}{-0.1cm}  
    \setlength{\belowcaptionskip}{-0.1cm} 
  \begin{center}  \includegraphics[width=0.85\textwidth]{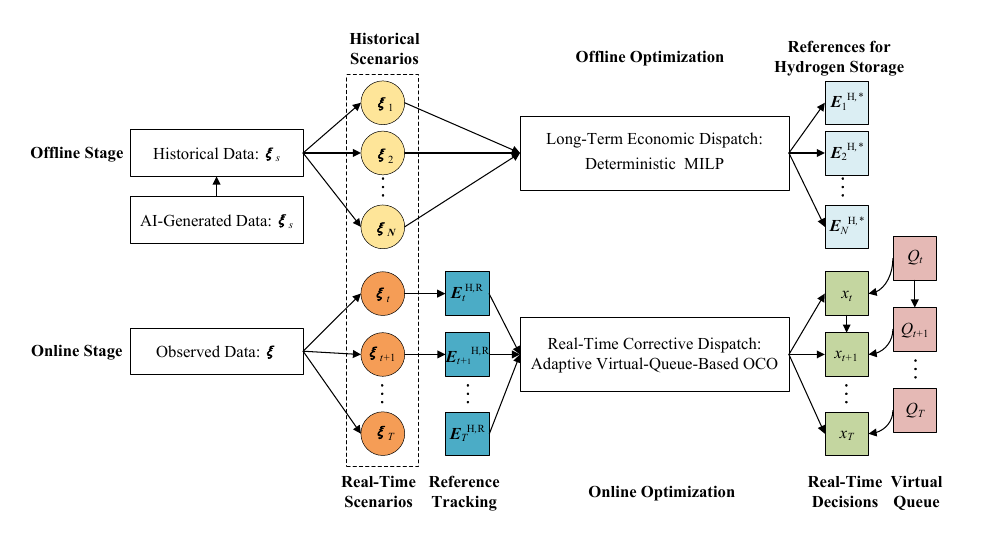}
     \caption{\rmfamily Diagram of prediction-free two-stage coordinated optimization framework.}\label{figureframework}
  \end{center}
  \vspace{-0.8cm}
\end{figure*}

We propose a two-stage coordinated optimization framework as illustrated in Fig~\ref{figureframework}. The proposed framework consists of both online and offline stage optimization. The offline stage aims to generate the ex-post SoC references for hydrogen storage using historical data on RES and load. These references can help avoid myopic decision-making and will be incrementally updated by kernel regression with newly observed data. Subsequently, online decisions are made using an adaptive virtual-queue-based OCO algorithm.
\subsection{Offline-stage optimization}
Firstly, sequential sequences of scenarios, denoted as $\bm{\xi_{s}}=\{\xi_{s,t}^{\text{L}},\xi_{s,t}^{\text{R}}\}$,$t\in{{\bm{\Omega} }_{T}}=\{1,2,\cdots,T\}, s\in{{\bm{\Omega} }_{S}}=\{1,2,\cdots,N\}$, are generated from historical data of the past few years. Additionally, to account for climate change and enhance the diversity of references, we can also collect references from different months and seasons. For instance, if we focus on a seasonal dispatch problem and have historical data for 5 years, then $T=1$ season and $N=5 \times 4$. To enhance adaptability to extreme weather conditions, we add extreme scenarios into the historical data using Generative Adversarial Networks~\cite{gui2021review}. Afterward, we can solve the deterministic 
mixed-integer linear programming (MILP) as~\eqref{P2} to generate the SoC references of hydrogen storage, i.e., $\bm{E}_{s}^{\text{H,*}}=\{E_{s,t}^{\text{H,*}}\}$,$t\in{{\bm{\Omega} }_{T}}, s\in{{\bm{\Omega} }_{S}}$.
\begin{alignat}{1}\label{P2}
  (\textbf{P}_{2})\hspace{4pt} \underset{\bm{x_{s}}}{\min} & {\hspace{4pt}} G(\bm{x_{s}},{\bm{\xi_{s}}})\notag\\ 
 \text{s.t.}  & {\hspace{4pt}} \eqref{hydrogen}\text{, }\eqref{DG}-\eqref{powerbalance}
\end{alignat}

\subsection{Online-stage optimization}

\textbf{(1) Data-Driven Reference Tracking}

Inspired by~\cite{guo2023long}, we propose a data-driven reference tracking method to combine both the 'lookback' pattern from historical data and the 'lookahead' pattern from newly observed data. Firstly, we define $\bm{\xi}_{[t]}$ as the observed sequence for uncertainties from the first time slot to the current time slot \textit{t} in~\eqref{observesequence}. Additionally, $\bm{\xi}_{s,[t]}$ defined in~\eqref{historicalsequence} represents the corresponding historical sequence for uncertainties in scenario \textit{s}. Subsequently, by checking the similarity between $\bm{\xi}_{[t]}$ and $\bm{\xi}_{s,[t]}$, dynamic weights $\omega_{s,t}$ are assigned to each historical scenario based on the Gaussian kernel function and Euclidean distance, as outlined in~\eqref{weight}. To account for the temporal dynamics, the Gaussian kernel function is modified with a scaling factor \textit{t}. And the optimal bandwidth \(\sigma\) can be found through heuristic methods such as the bisection method. Additionally, the weights are updated in real-time dispatch instead of using average or heuristic values. Finally, the SoC reference of hydrogen storage is updated as~\eqref{referweight}. This updated reference also determines the efficiency segment of hydrogen storage, eliminating the nonconvexity issue that arises when using convex optimization approaches.
\begin{subequations}\label{refupdate}
\begin{align}
    &\bm{\xi}_{[t]}=\{\xi_{1}^{\text{G}},\xi_{1}^{\text{L}},\xi_{1}^{\text{R}},\cdots,\xi_{t}^{\text{G}},\xi_{t}^{\text{L}},\xi_{t}^{\text{R}}\}
   \label{observesequence} \\
       &\bm{\xi}_{s,[t]}=\{\xi_{s,1}^{\text{G}},\xi_{s,1}^{\text{L}},\xi_{s,1}^{\text{R}},\cdots,\xi_{s,t}^{\text{G}},\xi_{s,t}^{\text{L}},\xi_{s,t}^{\text{R}}\}
   \label{historicalsequence} \\
   &\omega_{s,t}=\dfrac{K_{t}(\xi_{[t]},\xi_{s,[t]})}{\sum_{{s^{\prime}=1}}^{N}K_{t}(\xi_{[t]},\xi_{s^{\prime},{[t]}})},\ K_{t}(x,y)=e^{{-\frac{(\|x-y\|_{2})^{2}}{t\sigma^2}}}\label{weight}\\
      &\bm{E}_{[t]}^{\text{H,R}}=\sum\nolimits_{{s=1}}^{N}\omega_{s,t} \bm{E}_{s,[t]}^{\text{H,*}}\label{referweight}
\end{align}
\end{subequations}

\textbf{(2) Real-Time Corrective Dispatch}

Real-time corrective dispatch ($\textbf{P}_{3}$) is formulated in~\eqref{P3}, which aims to minimize the instant operational cost while tracking the SoC reference of hydrogen storage. $\varphi$ is the penalty coefficient to control the SoC deviation from the reference. $\textbf{P}_{3}$ admits a compact form in~\eqref{ft}. $f_t$ and $g_t$ represent the time-varying objective function and time-varying constraints due to hydrogen storage SoC reference ${E}_{t}^{\text{H,R}}$ and uncertainties $\bm{\xi}$, respectively. By leveraging the Lagrangian Relaxation, we can obtain the optimum by~\eqref{lanfun}. $\lambda_{t}$ is the dual variables of the constraints $g_t$. $\left\langle x,y\right\rangle$ denotes the standard inner product. However, without prior knowledge of uncertainties $\bm{\xi}$, $f_t$ and $g_t$ are unknown to the online decision-maker. Hence, we next design a VQB-OCO algorithm to solve this issue.
\begin{alignat}{1}\label{P3}
  (\textbf{P}_{3})\hspace{4pt} \underset{\bm{x}_{t}}{\min} & {\hspace{4pt}} G(\bm{x}_{t},{\bm{\xi}_{t}})+\varphi({E}_{t}^{\text{H}}-{E}_{t}^{\text{H,R}})^{2}\notag\\ 
 \text{s.t.}  & {\hspace{4pt}} \eqref{hydrogen}\text{, }\eqref{DG}-\eqref{powerbalance}
\end{alignat}
\vspace{-0.8cm}
\begin{subequations}
    \begin{align}
    & \min_{\bm{x}}\ f_t(x_t)\ \mathrm{~s.t.~}g_t(x_t)\leq0\label{ft}\\
    & x_t=\arg\min_{x}\{f_t(x)+\left\langle\lambda_t,\ g_t(x)\right\rangle \}\label{lanfun}
   \end{align}
\end{subequations}

\textbf{(3) VQB-OCO Algorithm}

The key idea of VQB-OCO is to use information from past time to approximate the current situation. The virtual queue is employed as the substitution of unknown dual variables. Hence, we design the update policy for the virtual queue, decisions, and weights in~\eqref{qt_update}-\eqref{weight_update}. Finally, we can obtain the weighted average value of dispatch decision as~\eqref{decision}. The VQB-OCO algorithm and the overall two-stage coordinated optimization framework are summarized in \textbf{Algorithm}~\ref{algorithm1}.

\begin{equation}\label{qt_update}
Q_{i,t-1}=Q_{i,t-2}+\beta_{t-1}[g_{t-1}(x_{t-1})]_+
\end{equation}
\begin{equation}\label{xt_update}
\begin{aligned}x_{i,t}&=\arg\min_{x}\{\alpha_{i,t-1}\left\langle\partial f_{t-1}(x_{t-1}),\ x\right\rangle+\\&\alpha_{i,t-1}\beta_{i,t-1}\left\langle Q_{i,t-1},\ [g_{t-1}(x)]_+\right\rangle+\|x-x_{i,t-1}\|^2\}\end{aligned}
\end{equation}
\begin{equation}\label{weight_update}
\ell_{i,t-1}=\left\langle\partial f_{t-1}(x_{t-1}),\ x_{i,t-1}-x_{t-1}\right\rangle,\ \rho_{i,t}=\dfrac{\rho_{i,t-1}e^{-\gamma \ell_{i,t-1}}}{\sum_{i=1}^N\rho_{i,t-1}e^{-\gamma \ell_{i,t-1}}}
\end{equation}
\begin{equation}\label{decision}
x_t=\sum\nolimits_{t=1}^N\rho_{i,t}x_{i,t}
\end{equation}

\begin{remark}[Approximation]
 $f_{t}(x)$ is approximated using the first-order Taylor expansion $\left\langle\partial f_{t-1}(x_{t-1}),\ x\right\rangle$. The term $\lambda_{t}$ is substituted by a virtual queue $Q_{i,t-1}$. The constraint function $g_t(x)$ is replaced by the clipped constraint function $\left[g_{t-1}(x)\right]_{+}$. A regularization term $\left|x-x_{i,t-1}\right|^2$ is added to ensure the convexity of the optimization problem and to enhance the convergence of the algorithm.
\end{remark}
\begin{remark}[Parallel Learning]
Determining the learning rate (step size) is important yet challenging. We assign different learning rates to the first two terms, $\alpha_{i,t-1}$ and $\beta_{i,t-1}$. Rather than utilizing fixed or adaptive learning rates, we employ the expert-tracking algorithm proposed by~\cite{zhang2018adaptive}, which computes $x_t$ in parallel with various learning rates as described in equation~\eqref{xt_update}. The weights for each expert $\rho_{i,t}$ are updated based on their empirical performance using an exponential function, as shown in equation~\eqref{weight_update}.

\end{remark}
\begin{remark}[Virtual Queue Updates]
Based on our previous work~\cite{qi2024chance}, the dual variables of the long-term constraints remain fixed when the optimum does not reach the constraint bounds. However, when the optimum reaches these bounds, the dual variables increase, representing a penalty. The update of the virtual queue follows the same pattern as described in equation~\eqref{qt_update} to limit constraint violations.
\end{remark}

\begin{algorithm}[htbp]\label{algorithm1}
\hspace{-0.4cm}\caption{Prediction-Free Two-Stage Online Optimization Algorithm}
\SetAlgoLined
\SetEndCharOfAlgoLine{}
\hspace{-0.4cm}\textbf{Stage1: Offline Optimization}

\hspace{-0.4cm}\KwIn{Historical scenarios of RES and load $\boldsymbol{\xi_{s}}$}
\hspace{-0.4cm}\KwOut{Historical reference for hydrogen storage $\bm{E}_{s}^{\text{H,*}}$}
\SetKwBlock{StepOne}{\hspace{-0.4cm}Step 1 -Initialization}{}
\SetKwBlock{StepTwo}{\hspace{-0.4cm}Step 2 - Reference Tracking \& VQB-OCO}{}
\SetKw{Parallel}{parallel}
    \For{$S=1$ \KwTo $N$ }{
        Solve the deterministic MILP problems ($\textbf{P}_{2}$) \; as~\eqref{P2} to generate the SoC references of\; hydrogen storage. \;
    }
\hspace{-0.4cm}\textbf{Stage2: Online Optimization}

\hspace{-0.4cm}\KwIn{Historical reference for hydrogen storage $\bm{E}_{s}^{\text{H,*}}$; \;  Real-time observation of RES and load $\boldsymbol{\xi_{t}}$}.
\hspace{-0.4cm}\KwOut{Real-Time Dispatch Decisions $x_{t}$}
\StepOne{
Set $Q_{i,t}=0,\ x_{i,1}\in \bm{X},\ x_1=\sum_{i=1}^N\rho_{i,1}x_{i,1}$,\\\vspace{0.1cm}$\rho_{i,1}=(M+1)/[i(i+1)M],\ \forall i\in\{1,2,\cdotp\cdotp\cdotp,M\}.$
}
\StepTwo{
    \For{$t=2$ \KwTo $T$ }{
      Update real-time SoC reference  as~\eqref{refupdate}; \;
      \For{$i=1$ \KwTo $M$\Parallel}{
      Update virtual queue $Q_{i,t}$ as~\eqref{qt_update};  \;
      Update decisions $x_{i,t}$ as~\eqref{xt_update};  \;
      Update weights $\rho_{i,t}$ as~\eqref{weight_update}.  \;
      }
      Calculate the dispatch decision $x_{i}$ as~\eqref{decision}.
    }
}
\end{algorithm}

\textbf{(4) Performance of VQB-OCO}

OCO focuses on the performance of \textit{regret} (Reg), as defined in~\eqref{reg}, where $y_t$ is the global optimum. Various OCO algorithms ensure that Reg is a sublinear function of \textit{T} by designing parameters and update policy, as it implies that the algorithm performs as well as the global optimum in hindsight as \textit{T} approaches infinity. Next, we provide parameter settings and a proof to achieve strictly sublinear dynamic regret.

\begin{equation}
\mathrm{Reg}=\sum\nolimits_{t=1}^T[f_t(x_t)-f_t(y_t)]\label{reg}
\end{equation}

\textit{Assumption 1.} The functions $f_{t}$ and $g_{t}$ are convex. The feasible set $\bm{X}$ is convex and closed, and it has a bounded diameter $d(\bm{X})$, i.e.,
\begin{equation}\label{assump1}
\parallel x-y\parallel\leq d(\bm{X}),\ \forall x,y\in \bm{X}
\end{equation}

\textit{Assumption 2.} There exists a positive constant 
\textit{F} such that
\begin{equation}\label{assump2}
\mid f_t(x)-f_t(y)\mid\leq F,\ \parallel g_t(x)\parallel\leq F,\ \forall t\in\bm{\Omega_{T}},\ \forall x,y\in \bm{X}
\end{equation}

\textit{Assumption 3.} The subgradients $\partial f_t(x)$ and $\partial g_t(x)$ exist. And there exists a positive constant \textit{G} such that
\begin{equation}\label{assump3}
\parallel\partial f_t(x)\parallel\leq G,\ \parallel\partial g_t(x)\parallel\leq G,\ \forall t\in\bm{\Omega_{T}},\ \forall x,y\in \bm{X}
\end{equation}

\begin{theorem}\label{performance}
Given the assumptions 1–3, and parameters setting as~\eqref{parameter}, $\kappa \in [0, c]$, $c \in (0, 1)$, $\alpha_0 > 0$, $\beta_0 > 0$, and $\gamma_0 \in (0, 1/(\sqrt{2G}))$ are constants. Then, we have the performance of Reg and Vio as~\eqref{results}. 
\end{theorem}
\begin{equation}\label{parameter}
\begin{aligned}
& M=\lfloor\kappa \log_{2}(1+T)\rfloor+1,\ \alpha_{i,t}=\dfrac{\alpha_{0}2^{i-1}}{t^{c}},\ \beta_{i,t}=\dfrac{\beta_{0}}{\sqrt{\alpha_{i,t}}}, \ \gamma=\dfrac{\gamma_{0}}{T^{c}}
\end{aligned}
\end{equation}
\begin{equation}\label{results}
\text{Reg}=\mathcal{O}(T^c (1 + P_x)^{1 - \kappa} + T^{1 - c} (1 + P_x)^{\kappa})
\end{equation}

\noindent\textit{\textbf{Proof:}} The performance of the proposed OCO algorithm achieves a similar performance with~\cite{muthirayan2022online} which achieves $\mathcal{O}(T^{\max \{1-a-c, c\}})$ for dynamic regret with the help of prediction data. And it outperforms the performance of~\cite{liu2022simultaneously} and \cite{ding2021dynamic}, achieving dynamic regret with a linear function of $P_x$, which is not satisfactory. Moreover, by setting $\kappa=c=0.5$, the proposed OCO algorithm achieve the performance of $\mathcal{O}(\sqrt{T(1 + P_x)})$, which aligns with the performance of~\cite{zhang2018adaptive} where
long term constraints are not considered. Hence, the proposed OCO algorithm is no worse than the existing versions. We defer the complete proof to \textbf{Appendix~\ref{proof1}}.

\section{Case studies}~\label{case}

\subsection{Set-up}

The main parameters and configurations are listed in Table~\ref{parameters}. Specifically, the capacities of the battery and hydrogen storage are half of the load capacity. The storage durations of the battery and hydrogen are 2 hours and 400 hours, respectively. The installed capacity of renewables is 200 kW, comprising an equal share of solar and wind. The cost coefficients can be found in~\cite{guo2023long}.

We demonstrate the effectiveness of the proposed method based on two datasets: (1) We use the 15-minute historical data on solar, wind, and load from 2014 to 2023 obtained from Belgium's transmission system operator (Elia)~\cite{Elia-data} for the baseline case study. (2) We also use the hourly historical data of wind and load from 1981 to 2020 in North China~\cite{North-China-data} to demonstrate the impact of data resolution and data quantity.

The optimization is coded in MatLab with Yalmip interface and solved by Gurobi 11.0 solver. The programming environment is Intel Core i9-13900HX @ 2.30GHz with RAM 32 GB.

\begin{table}[H]
\footnotesize\rmfamily
  \centering
  \caption{\rmfamily Parameters and configuration of the test microgrid.}
  \setlength{\tabcolsep}{0.8mm}{
      \begin{tabular}{c c c c c c}
    \toprule
    Parameters & Value & Parameters & Value & Parameters & Value \\
    \midrule
    Initial SoC  & 0.5   & $\overline{P}^{\text{H}}$  & 50 kW & $c^{\text{L}}$ & \$5/kWh \\
    $\eta^{\text{B,c/d}}$ & 0.9   & $\overline{E}^{\text{H}}$  & 20 MWh & $\overline{P}^{\text{D}}$ & 50kW \\
    $\epsilon$ & 1\%/month & $c^{\text{B}}$  & \$0.02/kWh & Wind Capacity & 100kW \\
    $\overline{P}^{\text{B}}$  & 50 kW &$c^{\text{H}}$  & \$0.03/kWh & Solar Capacity & 100kW \\
    $\overline{E}^{\text{B}}$   & 100kW & $c^{\text{D}}$ & \$0.3/kWh & Load Capacity & 100kW \\
    \bottomrule
    \end{tabular}%
    }\vspace{-0.3cm}\label{parameters}
\end{table}%
\FloatBarrier

\subsection{Offline-stage optimization}

(1) Data visualization

\begin{figure}[!ht]
\vspace{-0.5em}
 \footnotesize\rmfamily     \setlength{\abovecaptionskip}{-0.1cm}  
    \setlength{\belowcaptionskip}{-0.1cm} 
  \begin{center}  \includegraphics[width=0.9\columnwidth]{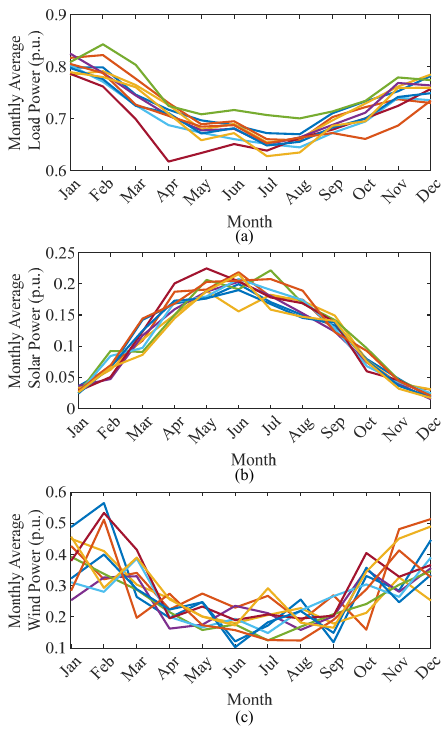}
     \caption{\rmfamily Data visualization on Elia dataset from 2014 to 2023: (a) monthly average load power, (b) monthly average solar power, and (c) monthly average wind power.}\label{datavisualization}
  \end{center}
  \vspace{-0.8cm}
\end{figure}

We first show the monthly average available renewable and load power of Elia from 2014 to 2023 in Figure~\ref{datavisualization}. It is observed that all of them exhibit seasonal patterns. Wind power is abundant in spring and winter but scarce in summer, while solar power is relatively high in summer and extremely low in winter. Load power peaks in winter. Correspondingly, the net load also peaks in winter and hits a low in summer. Therefore, it indicates the critical role of hydrogen storage to address the seasonal variations in renewables and load, as well as to maintain the long-term energy balance of the microgrid.

(2) Impact of hydrogen storage efficiency model

Next, we compare the offline energy management performance in 2023 with different hydrogen models, including:

\textbf{(E1):} Piecewise linear model as proposed in~\eqref{hydrogen}, and the parameters are fitted based on the experimental data as shown in Figure~\ref{chargeapproximate} and~\ref{dischargeapproximate}.

\textbf{(E2):} Constant efficiency model with both the highest charging and discharging efficiencies of 63\%.

\textbf{(E3):} Constant efficiency model with the lowest charging and discharging efficiencies, i.e., 53\% and 45\%, respectively.

\begin{figure}[!ht]
\vspace{-0.5em}
 \footnotesize\rmfamily     \setlength{\abovecaptionskip}{-0.1cm}  
    \setlength{\belowcaptionskip}{-0.1cm} 
  \begin{center}  \includegraphics[width=0.9\columnwidth]{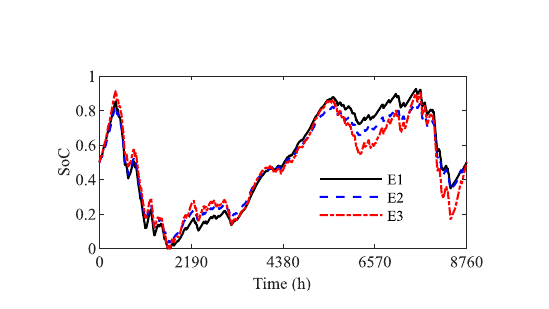}
     \caption{\rmfamily Hydrogen storage SoC over one representative year using different efficiency models.}\label{figefficiency}
  \end{center}
  \vspace{-0.8cm}
\end{figure}

The hydrogen storage SoC is shown in Figure~\ref{figefficiency}. Other results are also summarized in Table~\ref{tableefficiency}. It is observed that using the highest constant efficiency model results in the most optimistic performance, with the lowest operational cost (\$97,188), whereas the lowest constant efficiency model yields the highest operational cost (\$132,933) due to the significant increase in costs of diesel generation and loss of load. However, using a constant efficiency model can lead to feasibility issues in practical operation, resulting in losses in either charging or discharging power. Additionally, the optimistic strategy generated by the highest efficiency model will introduce an additional loss of load cost of 1.94 MWh. Considering the practical consequences, E2 and E3 will increase the total system costs by 10\% and 36\%, respectively, compared to E1. This result demonstrates that the proposed model can capture the characteristics of power-dependent efficiency and achieve more reliable and economical performance in practice.
\begin{table}[!ht]
\footnotesize\rmfamily
  \centering
  \caption{\rmfamily Yearly operational performance of the microgrid using different efficiency models.}
  \setlength{\tabcolsep}{0.8mm}{
      \begin{tabular}{c c c c c c c}
    \toprule
    \makecell{Model} & \makecell{$\text{Cost}^{\text{T/P}}$\\$(\$10^4)$} & \makecell{$\sum P_{t}^{\text{D,T/P}}\Delta t$\\$(\text{MWh})$} & \makecell{$\sum P_{t}^{\text{L,T/P}}\Delta t$\\$(\text{MWh})$} & \makecell{$\sum\Delta P_{t}^{\text{H,c}}\Delta t$\\$(\text{MWh})$} & \makecell{$\sum\Delta P_{t}^{\text{H,d}}\Delta t$\\$(\text{MWh})$} \\
        E1 & 9.87/9.87  & 324.68/324.68  & 0.00/0.00  & 0.00  & 0.00   \\
    E2 & 9.81/10.78  & 322.19/322.19  & 0.00/1.94  & 0.00  & -1.94      \\
    E3 & 13.29/13.29  & 374.33/374.19  & 3.85/3.85  & -0.14  & 0.00     \\
    \bottomrule
    \end{tabular}    }\vspace{-0.3cm}\label{tableefficiency}
\end{table}%

\subsection{Online-stage optimization}

(1) Reference tracking

We test the reference tracking performance in 2023 with different methods, including:

\textbf{(R1):} Global optimal reference generated by deterministic multi-period optimization with perfect knowledge of uncertainty realizations.

\textbf{(R2):} The proposed data-driven reference tracking, trained with 2014-2022 historical data.

\textbf{(R3):} The proposed data-driven reference tracking, trained with 2014-2022 historical data and AI-generated data. The AI-generated data is produced by randomly reducing the historical solar power but increasing the wind power by 10\%-50\% each quarter.

\textbf{(R4-R6):} Similar to the \textbf{(R3)} method, the AI-generated data is produced by:
- \textbf{(R4)}: Randomly reducing both the historical solar and wind power by 10\%-50\% each quarter.
- \textbf{(R5)}: Increasing the historical solar and wind power by 10\%-50\% each quarter.
- \textbf{(R6)}: Using all the AI-generated data from \textbf{(R3)} to \textbf{(R5)}.

\textbf{(R7):} The reference generated using the average historical performance.

\begin{figure}[!ht]
\vspace{-0.5em}
 \footnotesize\rmfamily     \setlength{\abovecaptionskip}{-0.1cm}  
    \setlength{\belowcaptionskip}{-0.1cm} 
  \begin{center}  \includegraphics[width=0.95\columnwidth]{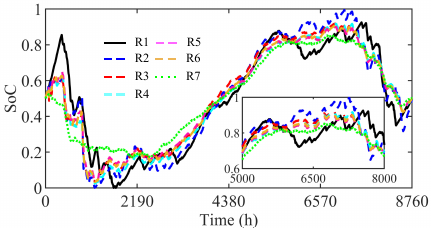}
     \caption{\rmfamily Hydrogen storage SoC references using different reference tracking methods.}\label{figreference}
  \end{center}
  \vspace{-0.8cm}
\end{figure}

The hydrogen SoC references are compared in Figure~\ref{figreference}, and the tracking performance is summarized in Table~\ref{tablereference}. The root mean square error (RMSE) is calculated as the average difference between the generated reference and the global optimal reference. The optimal choice of \(\sigma\) obtained through the bisection method is 0.098. It is observed that the references generated by the proposed methods \textbf{R2-R6} can better track the seasonal variations of RES and load, resulting in lower RMSE compared to the reference generated by \textbf{R7}. This is because the proposed methods employ kernel regression to update the weights of historical references instead of using fixed and average values. Additionally, additional generated data inputs will increase the tracking performance as they create new potential extreme scenarios. However, as is shown in the performance of \textbf{R6}, an excessive amount of generated data can lead to overfitting in the regression model, thereby reducing tracking accuracy. The average computation time for a single time interval is around 2 ms, which is acceptable even for minute-level scheduling and control.
\begin{table}[!ht]
\footnotesize\rmfamily
  \centering
  \caption{\rmfamily Reference tracking performance with different reference tracking methods.}
  \setlength{\tabcolsep}{1mm}{
    \begin{tabular}{cccccccc}
    \toprule
    Method & R1    & R2    & R3    & R4    & R5    & R6    & R7 \\
    \midrule
    RMSE  & ——    & 8.41\% & 8.96\% & 8.25\% & 8.01\% & 8.33\% & 10.54\% \\
    Time (ms) & ——    & 0.35  & 1.93  & 1.88  & 2.01  & 6.27  & 0.00  \\
    Data Size (Year) & ——    & 9     & 45    & 45    & 45    & 116   & 9 \\
    \bottomrule
    \end{tabular}    }\vspace{-0.3cm}\label{tablereference}
\end{table}%

Furthermore, we test the reference tracking performance on the North China dataset. The data visualization and reference tracking performance are shown in Figure~\ref{datavisualization1} in Appendix~\ref{A4}. Compared with the Elia dataset, both the wind and load data exhibit less variation across seasons and years. The maximum variations across years are 0.22 and 0.09 for wind and load, respectively, while for the Elia dataset, they are 0.25 and 0.15 for wind and load, respectively. Specifically, the load data in North China maintains the same shape across the years. The optimal choice of \(\sigma\) obtained through the bisection method is 50. The proposed tracking method performs the same as the averaged method, with an RMSE of 0.046. This is because the historical data and historical references show significant similarity across years, causing the proposed method to select the average value when updating weights. The above results demonstrate the benefit of using a data-driven reference tracking method when historical uncertainties exhibit significant variations across years. Additionally, an appropriate amount of AI-generated data can improve adaptability to extreme weather scenarios.

(2) Online decision-making

We test the online decision-making in 2023 with different dispatch methods, including:

\textbf{(M0):} Deterministic optimization with perfect knowledge of uncertainty realizations for the whole year results in the global optimum. This method is optimistic and not applicable in practice. Hence, it only serves as a baseline.

\textbf{(M1):} The proposed prediction-free coordinated approach, which utilizes OCO for online optimization and \textbf{R5} for reference tracking.

\textbf{(M2):} The scheduling-correction method proposed by~\cite{guo2023long}, which utilizes MPC for online optimization and \textbf{R5} for reference tracking.

\textbf{(M3):} Online optimization by OCO without reference tracking. 

\textbf{(M4):} Online optimization by MPC without reference tracking. 

\begin{figure}[!ht]
\vspace{-0.5em}
 \footnotesize\rmfamily     \setlength{\abovecaptionskip}{-0.1cm}  
    \setlength{\belowcaptionskip}{-0.1cm} 
  \begin{center}  \includegraphics[width=0.95\columnwidth]{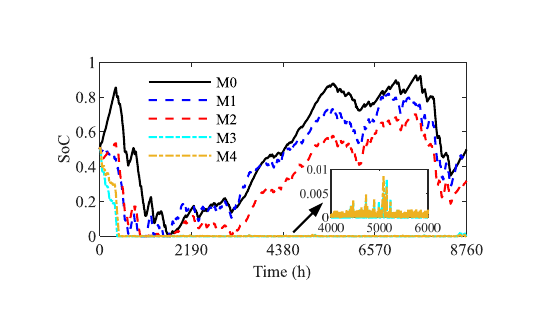}
     \caption{\rmfamily Hydrogen storage SoC strategies in Elia using different optimization methods.}\label{figreoptimization}
  \end{center}
  \vspace{-0.8cm}
\end{figure}

The operational performance is summarized in Table~\ref{tableoptimization}. It is observed that the proposed method \textbf{M1} outperforms the others in terms of cost-effectiveness, achieving an optimality gap of 27\% compared with \textbf{M0}. This is due to its smallest loss of load and RMSE compared to the reference. Additionally, \textbf{M1} and \textbf{M2} obtain much better economic performance than \textbf{M3} and \textbf{M4}. This can be explained by the hydrogen SoC as shown in Figure~\ref{figreoptimization}. \textbf{M3} and \textbf{M4} generate myopic decisions by continuously discharging the hydrogen storage to reduce short-term operational costs during the winter peak. They fail to charge the hydrogen storage during the renewable-rich spring and summer, resulting in an extremely low SoC after winter. Consequently, these myopic decisions prevent hydrogen storage from effectively shifting energy seasonally, leading to a substantial loss of load and low utilization of RES in practice. In contrast, \textbf{M1} and \textbf{M2} follow the pattern of reference while \textbf{M1} has the better reference following performance (lower RMSE) since OCO utilizes the real-time observed data. This result demonstrates the benefit of introducing a global reference for the online optimization method. We also apply SDP algorithm
in~\cite{darivianakis2017data} to this problem. However, the SDP cannot converge within 24 hr due to the ``curse of dimensionality''. Therefore, it is infeasible to use SDP for long-term energy management of microgrid with H-BES.

\begin{figure}[!ht]
\vspace{-0.5em}
 \footnotesize\rmfamily     \setlength{\abovecaptionskip}{-0.1cm}  
    \setlength{\belowcaptionskip}{-0.1cm} 
  \begin{center}  \includegraphics[width=0.95\columnwidth]{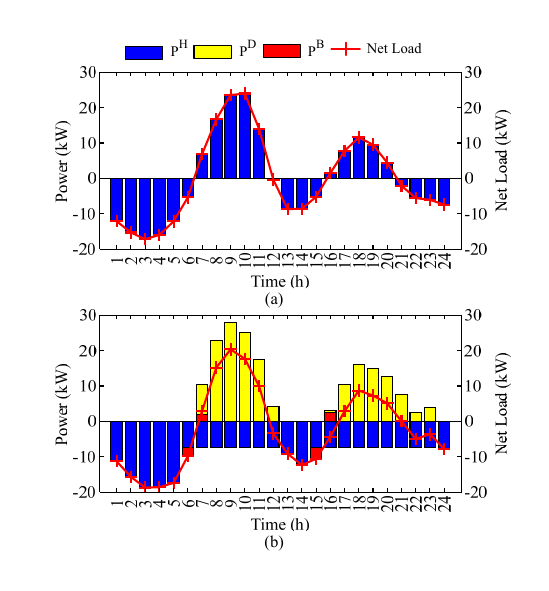}
     \caption{\rmfamily Comparison of power dispatch strategies with different optimization methods: (a) M1 and (b) M2.}\label{figrecomposition}
  \end{center}
  \vspace{-0.8cm}
\end{figure}

Moreover, we compare the power dispatch strategies of H-BES and DG using \textbf{M1} and \textbf{M2}, as shown in Figure~\ref{figrecomposition}. It is observed that \textbf{M1} can better track the net load curve using only hydrogen storage actions. In contrast, \textbf{M2} keeps charging hydrogen storage and uses DG when renewables are insufficient. This is because the OCO-based method simultaneously tracks the previous decisions and the reference, updating the strategy based on newly observed data, which is more adaptive to the time-varying environment. While the MPC-based method only tracks the reference and updates the strategy based on forecast data. Therefore, if the reference or forecast is not accurate, the MPC-based method may struggle to achieve good performance. This can also be explained by SoC gaps as shown in Figure~\ref{figreoptimization}. Furthermore, due to prediction errors, MPC-based online optimization may encounter infeasibility issues, resulting in additional loss of load and penalty costs. In contrast, the OCO-based method makes decisions based on observed data, thereby avoiding infeasibility issues. Regarding computational efficiency, it can be seen that the OCO method has better performance than MPC. Both methods achieve single-step optimization in tens of ms, which is acceptable for most online optimization scenarios.

\begin{table}[!ht]
\footnotesize\rmfamily
  \centering
  \caption{\rmfamily Yearly operational performance of the microgrid in Elia using different optimization methods.}
  \setlength{\tabcolsep}{0.3mm}{
      \begin{tabular}{c c c c c c c}
    \toprule
    \makecell{Method} & \makecell{$\text{Cost}$ $(\$10^4)$} & \makecell{$\sum P_{t}^{\text{D}}\Delta t$ $(\text{MWh})$} & \makecell{$\sum P_{t}^{\text{L}}\Delta t$ $(\text{MWh})$} & RMSE & \makecell{ Time (ms)}\\
    M0    & 9.87 & 324.68 & 0.00 & 0.00  & 25.20  \\
    M1    & 12.55 & 336.30 & 4.59 & 10.46  & 87.87  \\
    M2    & 17.38 & 488.46 & 5.19 & 19.70  & 97.77  \\
    M3    & 29.43 & 257.77 & 43.19 & 53.19  & 48.59  \\
    M4    & 38.05 & 457.58 & 48.45 & 52.93  & 51.64  \\
    \bottomrule
    \end{tabular}    }\vspace{-0.3cm}\label{tableoptimization}
\end{table}%

Additional tests on the North China dataset align with the above results, showing that optimization with reference outperforms optimization without reference, and OCO-based methods outperform MPC-based methods. The results are summarized in Table~\ref{tableoptimizationChina} and Figure~\ref{figreoptimizationChina} in Appendix~\ref{A4}. Compared with the case results in Elia, both \textbf{M1} and \textbf{M2} achieve SoC strategies with smaller gaps from \textbf{M0}, due to better reference tracking performance. However, in terms of cost and reliability performance, the North China case shows worse results. This is because renewable energy in the North China case is solely supplied by wind power, leading to insufficient generation and higher load curtailment.

\subsection{Sensitivity analysis}

In this subsection, we further investigate the key impact factor of the proposed optimization framework.

\begin{figure}[!ht]
\vspace{-0.5em}
 \footnotesize\rmfamily     \setlength{\abovecaptionskip}{-0.1cm}  
    \setlength{\belowcaptionskip}{-0.1cm} 
  \begin{center}  \includegraphics[width=0.95\columnwidth]{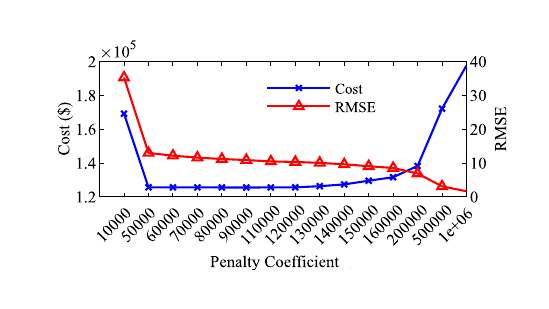}
     \caption{\rmfamily Reliability performance of microgrid with different renewable capacity using M1 and M2 methods.}\label{figurepenalty}
  \end{center}
  \vspace{-0.8cm}
\end{figure}

\textbf{(1) Penalty Coefficient of Reference Tracking.} The penalty coefficient represents the tradeoff between instant operational cost and reference tracking performance. However, since this reference is estimated from historical data, it may not be optimal for the current year. We compare the cost and tracking performance in Figure~\ref{figurepenalty} when scaling up the penalty coefficient $\varphi$. It is observed that RMSE is monotonically decreased with the penalty coefficient, while the operational cost initially decreases to a minimum value at \(\varphi = 90000\) but then gradually increases with the penalty coefficient. This suggests that improving tracking performance does not always lead to lower operational costs.

\textbf{(2) Reference and step size of OCO Algorithm.} As illustrated in Section 5.3, the reference has a critical impact on operational performance. Table~\ref{tableOCO} summarizes the performance using M1 with different references (i.e., fixed reference generated by \textbf{R7} and updated reference generated by \textbf{R5}) and different step sizes (i.e., fixed step size proposed by~\cite{muthirayan2022online} and step size generated by the proposed expert-tracking algorithm). It is observed that compared to fixed reference, using the proposed updated reference reduces operational cost by 1.09\%-1.70\% and the loss of load by 5.68\%-8.52\%. This highlights the importance of finding the right reference for hydrogen storage. Additionally, the proposed step size setting decreases the operational cost by 0.67\%-1.29\% and the loss of load by 0.06\%-3.16\%. This is because the step size generated by expert-tracking can better adapt to the changing cost function and avoid heuristic settings.

\begin{table}[!ht]
\footnotesize\rmfamily
  \centering
  \caption{\rmfamily Yearly operational performance of the microgrid in Elia using M1 with different references and step sizes.}
  \setlength{\tabcolsep}{0.8mm}{
      \begin{tabular}{c c c c c c c}
    \toprule
    \makecell{Reference} & \makecell{step size} &\makecell{$\text{Cost}^{\text{T/P}}$\\$(\$10^4)$} & \makecell{$\sum P_{t}^{\text{D,T/P}}\Delta t$\\$(\text{MWh})$} & \makecell{$\sum P_{t}^{\text{L,T/P}}\Delta t$\\$(\text{MWh})$} & RMSE \\
\toprule
    \multirow{2}[2]{*}{Updated} & Fixed & 12.55  & 336.30  & 4.60  & 10.87  \\
          & Expert-Tracking & 12.47  & 333.66  & 4.60  & 10.99  \\
    \midrule
    \multirow{2}[2]{*}{Fixed} & Fixed & 12.77  & 336.33  & 5.03  & 9.26  \\
          & Expert-Tracking & 12.61  & 333.58  & 4.87  & 9.26  \\
    \bottomrule
    \end{tabular}    }\vspace{-0.3cm}\label{tableOCO}
\end{table}%

\begin{figure}[!ht]
\vspace{-0.5em}
 \footnotesize\rmfamily     \setlength{\abovecaptionskip}{-0.1cm}  
    \setlength{\belowcaptionskip}{-0.1cm} 
  \begin{center}  \includegraphics[width=0.95\columnwidth]{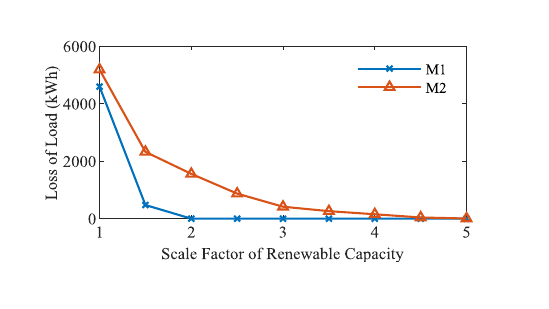}
     \caption{\rmfamily Reliability performance of microgrid with different renewable capacity using M1 and M2 methods.}\label{figurescale}
  \end{center}
  \vspace{-0.8cm}
\end{figure}

\textbf{(3) Sizing of Renewables.} We further compare the reliability performance of \textbf{M1} and \textbf{M2} in Figure~\ref{figurescale} when scaling up renewable capacity. It is observed that the loss of load decreases with increased renewable capacity. Moreover, \textbf{M1} achieves an acceptable reliability level with twice the renewable capacity, while \textbf{M2} requires at least 4 times the renewable capacity to meet reliability requirements. This demonstrates the benefit of the proposed method in reducing renewable planning costs. 

\section{Conclusion}~\label{conclusion}

This paper proposes a prediction-free coordinated optimization framework for long-term energy management of microgrid with H-BES. To accurately captures the power-dependent efficiency of hydrogen storage, we propose an approximate semi-empirical hydrogen storage model using piecewise linear relaxation. Moreover, to address the long-term operational patterns of renewables and load and to eliminate dependence on predictions, we introduce a prediction-free, two-stage coordinated optimization framework. The key idea is to generate and track the SoC reference of hydrogen storage using historical scenarios and kernel regression to avoid myopic online decisions. And online decisions are made based on the proposed VQB-OCO algorithm by leveraging feedback control policy and newly observed data. Case studies on Elia and North China verify that:

(1) Compared to the constant efficiency model, the proposed approximation model avoids both overly optimistic and overly conservative strategies.

(2) The proposed optimization framework outperforms existing online optimization methods and achieves an acceptable gap compared to deterministic optimization with perfect foresight of uncertainties.

(3) The SoC reference of hydrogen storage is critical to overall performance. Therefore, more reliable historical or AI-generated scenarios, along with sophisticated techniques for setting penalty coefficients, are highly required.

(4) The OCO algorithm typically lacks a global view and is sensitive to step size settings. Thus, it is beneficial to incorporate a penalty term for long-term pattern tracking and expert-tracking for step size updates.

Further work will focus on addressing constraints violations in the OCO algorithm and extending the proposed framework to operations in a market environment.








\begin{center}
     \textbf{Appendix}
\end{center}

\begin{appendices}

\section{Proof of storage priority}\label{priority}

The marginal discharge cost of battery is \$10/MWh-\$30/MWh, while it is \$1/kg-\$2/kg (\$30/MWh-\$60/MWh) for hydrogen storage, which is much higher than battery. Assume we have an optimal discharge power from H-BES, which is entirely supplied by the battery, i.e., $x^\text{1}=P_t^{\text{B,d}}$. Considering another optimum where the discharge power is a mix from both the battery and hydrogen storage, i.e., $x^\text{2}=\rho P_t^{\text{B,d}}+(1-\rho) P_t^{\text{H,d}}$. From the H-BES cost of the two optima~\eqref{storage priority} and the marginal cost \& efficiency of two types of storages, we can draw the conclusion that $C_t^\text{H-BES,2}>C_t^\text{H-BES 1}$. This indicates that hydrogen storage will not be activated until the battery is fully discharged or charged. Hence, we finish the proof.
\begin{subequations}\label{storage priority}
\begin{align}
    &\hspace{-0.7cm}P_t^{\text{B,d}}=\rho P_t^{\text{B,d}}+(1-\rho) P_t^{\text{H,d}}   \label{pe} \\
      &\hspace{-0.7cm}C_t^\text{H-BES,1}=c^\text{B}P_{t}^\text{B,d}\Delta t\\
      &\hspace{-0.7cm}C_t^\text{H-BES,2}=(c^\text{B}\rho P_{t}^\text{B,d}+c^\text{H}(1-\rho)P_{t}^\text{H,d})\Delta t \label{costcompare}
\end{align}
\end{subequations}

\section{Proof of bounded dynamic regret}\label{proof1}
Let $\{x_{i,t}\}$ and $\{x_{t}\}$ be the sequences generated by Algorithm~\ref{algorithm1}. Let $\{y_t\}$ be a global optimum in the feasible set $\bm{X}$. From $f_{t}$ is convex and~\eqref{assump3}, we have:
\begin{equation}\label{reg1}
\begin{aligned}
&\hspace{-0.3cm}f_t(x_{i,t}) - f_t(y_t) \leq \left\langle \partial f_t(x_{i,t}),\ x_{i,t} - y_t \right\rangle \\
&\hspace{-0.3cm}\leq G \parallel x_{i,t} - x_{i,t+1} \parallel + \left\langle \partial f_t(x_{i,t}),\ x_{i,t+1} - y_t \right\rangle \\
&\hspace{-0.3cm}\leq \frac{G^{2} \alpha_{i,t}}{2} + \frac{1}{2 \alpha_{i,t}} \parallel x_{i,t} - x_{i,t+1} \parallel^{2} + \left\langle \partial f_t(x_{i,t}),\ x_{i,t+1} - y_t \right\rangle
\end{aligned}
\end{equation}
For the rightmost term of~\eqref{reg1}, we have:
\begin{equation}\label{reg2}
\begin{aligned}
\begin{aligned}&\left\langle\partial f_t(x_{i,t}),\ x_{i,t+1}-y_t\right\rangle 
\\&=\left\langle\beta_{i,t+1}(\partial[g_t(x_{i,t+1})]_{+})^{T}Q_{i,t},\ y_t-x_{i,t+1}\right\rangle 
\\&+\left\langle\partial f_t(x_{i,t})+\beta_{i,t+1}(\partial[g_t(x_{i,t+1})]_{+})^{T}Q_{i,t},\ x_{i,t+1}-y_t\right\rangle\end{aligned}
\end{aligned}
\end{equation}
Since $g_{t}$ is a convex function, it is trivial to show that $[g_t]_+$ is also convex; hence the first term of~\eqref{reg2} can be relaxed:
\begin{equation}\label{reg3}
\begin{aligned}&\left\langle\beta_{t+1}(\partial[g_t(x_{i,t+1})]_{+})^{T}Q_{i,t},\ y_t-x_{i,t+1}\right\rangle\\&\leq\beta_{t+1}\left\langle Q_{i,t},\ [g_t(y_t)]_+\right\rangle-\beta_{i,t+1}\left\langle Q_i(t),\ [g_t(x_{i,t+1})]_+\right\rangle\end{aligned}
\end{equation}
From Lemma 1 in~\cite{yi2020distributed}, we have:
\begin{equation}\label{reg4}
\begin{aligned}
&\left\langle\partial f_t(x_{i,t})+\beta_{i,t+1}(\partial[g_t(x_{i,t+1})]_{+})^{T}Q_{i,t},\ x_{i,t+1}-y_t\right\rangle\\&\leq\frac1{\alpha_{i,t}}(\parallel y_t-x_{i,t}\parallel^2-\parallel y_t-x_{i,t+1}\parallel^2-\parallel x_{i,t+1}-x_{i,t}\parallel^2)
\end{aligned}
\end{equation}
Combining~\eqref{weight_update},~\eqref{reg1}-\eqref{reg4}, we have:
\begin{equation}\label{reg5}
\begin{aligned}
\hspace{-0.4cm}\ell_t(x_{i,t}) - \ell_t(y_t)& \leq\frac{G^2\alpha_{i,t}}2+\frac1\alpha(\parallel y_t-x_{i,t}\parallel^2-\parallel y_t-x_{i,t+1}\parallel^2)
\\&+\beta_{t+1}\left\langle Q_i(t),\ [g_t(y_t)]_+\right\rangle
\end{aligned}
\end{equation}
Since the last term of~\eqref{reg5} is non-negative, we have:
\begin{equation}\label{reg6}
\begin{aligned}
&\sum_{t=1}^{T} (\ell_t(x_{i,t}) - \ell_t(y_t))\leq\sum_{t=1}^T\frac{G^2\alpha_{i,t}}2\\
&+\sum_{t=1}^T\frac1{\alpha_{i,t}}(\parallel y_t-x_{i,t}\parallel^2-\parallel y_t-x_{i,t+1}\parallel^2)
\end{aligned}
\end{equation}
For the first term of~\eqref{reg6}, we have:
\begin{equation}\label{reg7}
\begin{aligned}
\sum_{t=1}^T\frac{G^2\alpha_{i,t}}2\leq\frac{2^{i-1}G^2}2\sum_{t=1}^T\frac1{t^{c}}\leq\frac{2^{i-1}G^2}{2(1-c)}T^{1-c}
\end{aligned}
\end{equation}
By leveraging~\eqref{assump1} and update policy~\eqref{parameter}, we have:
\begin{equation}\label{reg8}
\begin{aligned}
&\sum_{t=1}^{T} \frac{t^c}{\alpha_0 2^{i-1}} \left( \|y_t - x_{i,t}\|^2 - \|y_t - x_{i,t+1}\|^2 \right) \\
&= \frac{1}{\alpha_0 2^{i-1}}\sum_{t=1}^{T} \big( t^c \|y_t - x_{i,t}\|^2 - (t + 1)^c \|y_{t+1} - x_{i,t+1}\|^2 \\
&+  (t + 1)^c \|y_{t+1} - x_{i,t+1}\|^2 - t^c \|y_t - x_{i,t+1}\|^2  \\
&+ t^c \|y_t - x_{i,t+1}\|^2 - t^c \|y_t - x_{i,t}\|^2 \big)\\
&\le \frac{1}{\alpha_0 2^{i-1}} \|y_1 - x_{i,1}\|^2 + \frac{1}{\alpha_0 2^{i-1}} \sum_{t=1}^{T} \left( (t + 1)^c - t^c \right) (d(\mathbb{X}))^2 \\
&+ \frac{2}{\alpha_0 2^{i-1}} \sum_{t=1}^{T} t^c d(\bm{X}) \|y_{t+1} - y_t\| \\
&\le \frac{1}{\alpha_0 2^{i-1}} \left( 1 + (T + 1)^c - 1 \right) (d(\bm{X}))^2 + \frac{2T^c d(\bm{X}) P_x}{\alpha_0 2^{i-1}} \\
&\le \frac{2}{\alpha_0 2^{i-1}} (d(\bm{X}))^2 T^c \left( 1 + \frac{P_x}{d(\bm{X})} \right)
\end{aligned}
\end{equation}
Let $i_0=\left\lfloor\frac12\log_2(1+\frac{P_x}{d(\bm{X})})\right\rfloor+1\in[N]$, such that we have:
\begin{equation}\label{reg9}
\begin{aligned}
2^{i_0-1}\leq\sqrt{1+\frac{P_x}{d(\bm{X})}}\leq2^{i_0}.
\end{aligned}
\end{equation}
Combining~\eqref{reg7}-\eqref{reg9} yields:
\begin{equation}\label{reg10}
\begin{aligned}
\sum_{t=1}^{T} (\ell_t(x_{i_0, t}) - \ell_t(y_t))& \leq \frac{4}{\alpha_0} (d(\bm{X}))^2 T^c \left( 1 + \frac{P_x}{d(\bm{X})} \right)^{1-\kappa}\\
&+ \frac{G^2 \alpha_0}{2 (1 - c)} T^{1-c} \left( 1 + \frac{P_x}{d(\bm{X})} \right)^{\kappa}
\end{aligned}
\end{equation}
Applying Lemma 1 in reference~\cite{zhang2018adaptive} to~\eqref{weight_update} and~\eqref{decision} yields:
\begin{equation}\label{reg11}
\begin{aligned}
\sum_{t=1}^T\ell_t(x_t)-\min_{i\in[N]}\{\sum_{t=1}^T\ell_t(x_{i,t})+\frac1\gamma\ln\frac1{\rho_{i,1}}\}\leq\frac{\gamma(Gd(\bm{X}))^2T}2
\end{aligned}
\end{equation}
\begin{equation}\label{reg12}
\begin{aligned}
\sum_{t=1}^T(\ell_t(x_t)-\ell_t(x_{i_0,t}))\leq\frac{\gamma_0 (G d(\bm{X}))^2 T^{1 - c}}{2} + \frac{1}{\gamma_0} T^c \ln \frac{1}{\rho_{i_0, 1}}
\end{aligned}
\end{equation}
From $\rho_{i,1}=(M+1)/[i(i+1)M]$, we have:
\begin{equation}\label{reg13}
\begin{aligned}
\ln\frac1{\rho_{i_0,1}}\leq\ln(i_0(i_0+1))\leq2\ln(i_0+1)\leq2\ln(\left\lfloor\kappa\log_2(1+\frac{P_x}{d(X)})\right\rfloor)
\end{aligned}
\end{equation}
From~\eqref{weight_update} and that $f_{t}$ is convex, we have
\begin{equation}\label{reg14}
\begin{aligned}
f_t(x_t)-f_t(y_{t})\leq \ell_t(x_t)-\ell_t(y_{t})
\end{aligned}
\end{equation}
Combining~\eqref{reg10}-\eqref{reg14} yields:
\begin{equation}\label{reg15}
\begin{aligned}
&\text{Reg}\leq \frac{4}{\alpha_0} (d(\bm{X}))^2 T^c \left( 1 + \frac{P_x}{d(\bm{X})} \right)^{1 - \kappa}+\frac{\gamma_0 (G d(\bm{X}))^2 T^{1 - c}}{2}\\
 &+\frac{G^2 \alpha_0}{2(1 - c)} T^{1 - c} \left( 1 + \frac{P_x}{d(\bm{X})} \right)^{\kappa}+\frac{2}{\gamma_0} T^c \ln ([\kappa \log_2 \left( 1 + \frac{P_x}{d(\bm{X})} \right)])\\
\end{aligned}
\end{equation}
Hence, we finish the proof.

\section{Results on North China Dataset}\label{A4}

\begin{figure}[!ht]
\vspace{-0.5em}
 \footnotesize\rmfamily     \setlength{\abovecaptionskip}{-0.1cm}  
    \setlength{\belowcaptionskip}{-0.1cm} 
  \begin{center}  \includegraphics[width=0.9\columnwidth]{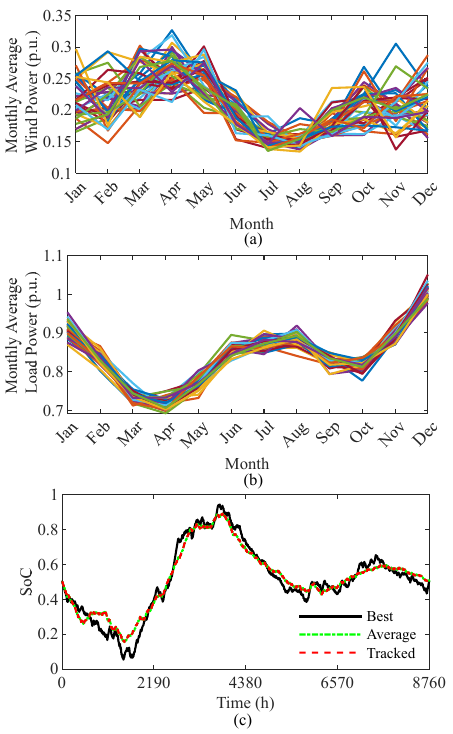}
     \caption{\rmfamily Data visualization and reference tracking on North China dataset: (a) monthly average wind power, (b) monthly average load power, and (c) hydrogen storage SoC reference.}\label{datavisualization1}
  \end{center}
  \vspace{-0.8cm}
\end{figure}

\begin{table}[!ht]
\footnotesize\rmfamily
  \centering
  \caption{\rmfamily Yearly operational performance of the microgrid in North China using different optimization methods.}
  \setlength{\tabcolsep}{0.5mm}{
      \begin{tabular}{c c c c c c c}
    \toprule
    \makecell{Method} & \makecell{$\text{Cost}$ $(\$10^5)$} & \makecell{$\sum P_{t}^{\text{D}}\Delta t$ $(\text{MWh})$} & \makecell{$\sum P_{t}^{\text{L}}\Delta t$ $(\text{MWh})$} & RMSE & \makecell{ Time (ms)}\\
    M0    & 5.30  & 411.78  & 80.60  & 0.00  & 24.28  \\
    M1    & 11.74  & 424.92  & 208.85  & 0.08  & 78.28  \\
    M2    & 15.49  & 257.52  & 1289.70  & 0.10  & 85.39  \\
    M3    & 13.35  & 349.87  & 245.83  & 53.19  & 35.46  \\
    M4    & 23.24  & 217.91  & 2104.74  & 52.93  & 58.31  \\
    \bottomrule
    \end{tabular}    }\vspace{-0.3cm}\label{tableoptimizationChina}
\end{table}%

\begin{figure}[!ht]
\vspace{-0.5em}
 \footnotesize\rmfamily     \setlength{\abovecaptionskip}{-0.1cm}  
    \setlength{\belowcaptionskip}{-0.1cm} 
  \begin{center}  \includegraphics[width=0.95\columnwidth]{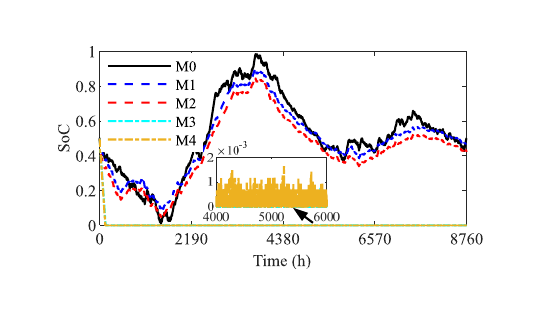}
     \caption{\rmfamily Hydrogen storage SoC strategies in North China using different optimization methods.}\label{figreoptimizationChina}
  \end{center}
  \vspace{-0.8cm}
\end{figure}

\end{appendices}

\printcredits
\vspace{0.3cm}

\noindent \textbf{Declaration of competing interest}

The authors declare that they have no known competing financial interests or personal relationships that could have appeared to influence the work reported in this paper.
\vspace{0.3cm}

\noindent \textbf{Data availability}

The original data can be downloaded from~\cite{Elia-data} and~\cite{North-China-data}. 
\vspace{0.3cm}

\noindent \textbf{Acknowledgements}

This work is partly supported by the National Science Foundation under award ECCS-2239046 and partly supported by the Columbia University SEAS Interdisciplinary Research Seed (SIRS) fund. Ning Qi graciously acknowledges special funding from the China Postdoctoral Science Foundation (No.2023TQ0169).

\bibliographystyle{LDES}

\bibliography{LDES}



\end{document}